\documentclass{amsart}

\usepackage{graphicx}
\usepackage{subfigure}
\usepackage{amsmath}
\usepackage{amsthm}
\usepackage{amsfonts}

\usepackage{hyperref}

\usepackage{verbatim}

\theoremstyle{definition}

\theoremstyle{remark}

\numberwithin{equation}{section}

\begin{document}

\title{Spectrum of the Laplacian on Regular Polyhedra}

\author{Evan Greif}
\address{Department of Mathematics, Harvard University, 1 Oxford St., Cambridge, MA, USA}
\email{evangreif@college.harvard.edu}
\thanks{EG was supported by the National Science Foundation through the Research Experience for Undergraduates (REU) Program, Grant DMS-1156350.}

\author{Daniel Kaplan}
\address{Department of Mathematics, Imperial College London, 180 Queens Gate, Kensington, London SW7 2RH, United Kingdom}
\email{d.kaplan15@imperial.ac.uk}
\thanks{DK was supported by the National Science Foundation through the Research Experience for Undergraduates (REU) Program, Grant DMS-1156350.}

\author{Robert S. Strichartz}
\address{Department of Mathematics, Cornell University, Malott Hall, Ithaca, NY 14853, USA}
\email{str@math.cornell.edu}
\thanks{RSS was supported in part by the National Science Foundation, Grant DMS-1162045.}

\author{Samuel C. Wiese}
\address{Department of Mathematics, Universit\"at Leipzig, Augustusplatz 10, 04109 Leipzig, Germany}
\email{sw31hiqa@studserv.uni-leipzig.de}
\thanks{SCW was supported by the Foundation of German Business (SDW)}


\date{August 13, 2018.}


\begin{abstract}

We study eigenvalues and eigenfunctions of the Laplacian on the surfaces of four of the regular polyhedrons: tetrahedron, octahedron, icosahedron and cube. We show two types of eigenfunctions: nonsingular ones that are smooth at vertices, lift to periodic functions on the plane and are expressible in terms of trigonometric polynomials; and singular ones that have none of these properties. We give numerical evidence for conjectured asymptotic estimates of the eigenvalue counting function. We describe an enlargement phenomenon for certain eigenfunctions on the octahedron that scales down eigenvalues by a factor of $\frac{1}{3}$.

\end{abstract}

\maketitle


\section{Introduction.}

The regular polyhedra are the two-dimensional surfaces that bound the Platonic solids. This is such an ancient topic in mathematics that one might think that there is nothing new to be said about it. Perhaps this is so when it comes to questions of geometry. But the purpose of this paper is to highlight new and interesting questions on analysis, specifically concerning eigenvalues and eigenfunctions of the Laplacian on regular polyhedra. For example, Figures 1a and 1b show images of two eigenfunctions on the cube. The first we call a \emph{nonsingular eigenfunction}. It extends to a smooth periodic function on the plane, and it has a formula given by trigonometric polynomials. The second we will call a \emph{singular eigenfunction}. It is not smooth at the vertices of the cube (see Figure 2a and 2b for the graph of the restriction to line segments passing through vertices). We believe that these functions should belong to a new class of special functions. For another interesting example of \emph{nonsingular} eigenfunctions on the octahedron see Figures 3a and 3b. The eigenfunction in Figure 3a is obtained from the one in Figure 3b by rotating and dilating by the factor $\sqrt{3}$, which reduces the eigenvalue by a factor of $\frac{1}{3}$.

This paper is addressed to a general audience, so we will sidestep some technical issues. Our surfaces are represented by regions in the plane with edges identified. See Figure $4$ for the cube.

The Laplacian on the surface is just the usual two dimensional Laplacian $\Delta = \frac{\partial^2}{\partial x^2} + \frac{\partial^2}{\partial y^2}$ on the planar realization, such that values and normal derivatives along the identified edges match. At the singular vertices we don't impose a differential equation; it is enough to assume that the function is continuous. This is a special case of a general theory of Laplacians on Alexandrov spaces. (For the experts please see [4].) By the spectrum we mean a study of both the eigenvalues $\lambda$ and and eigenfunctions $u$ satisfying 
\begin{equation}
\label{eq:cond}
-\Delta u = \lambda u\text{.}
\end{equation}


\begin{figure}[h]
    \centering
    \subfigure[For eigenvalue $\lambda_{37}=8$ (nonsingular)]
    {
        \includegraphics[width=2.2in]{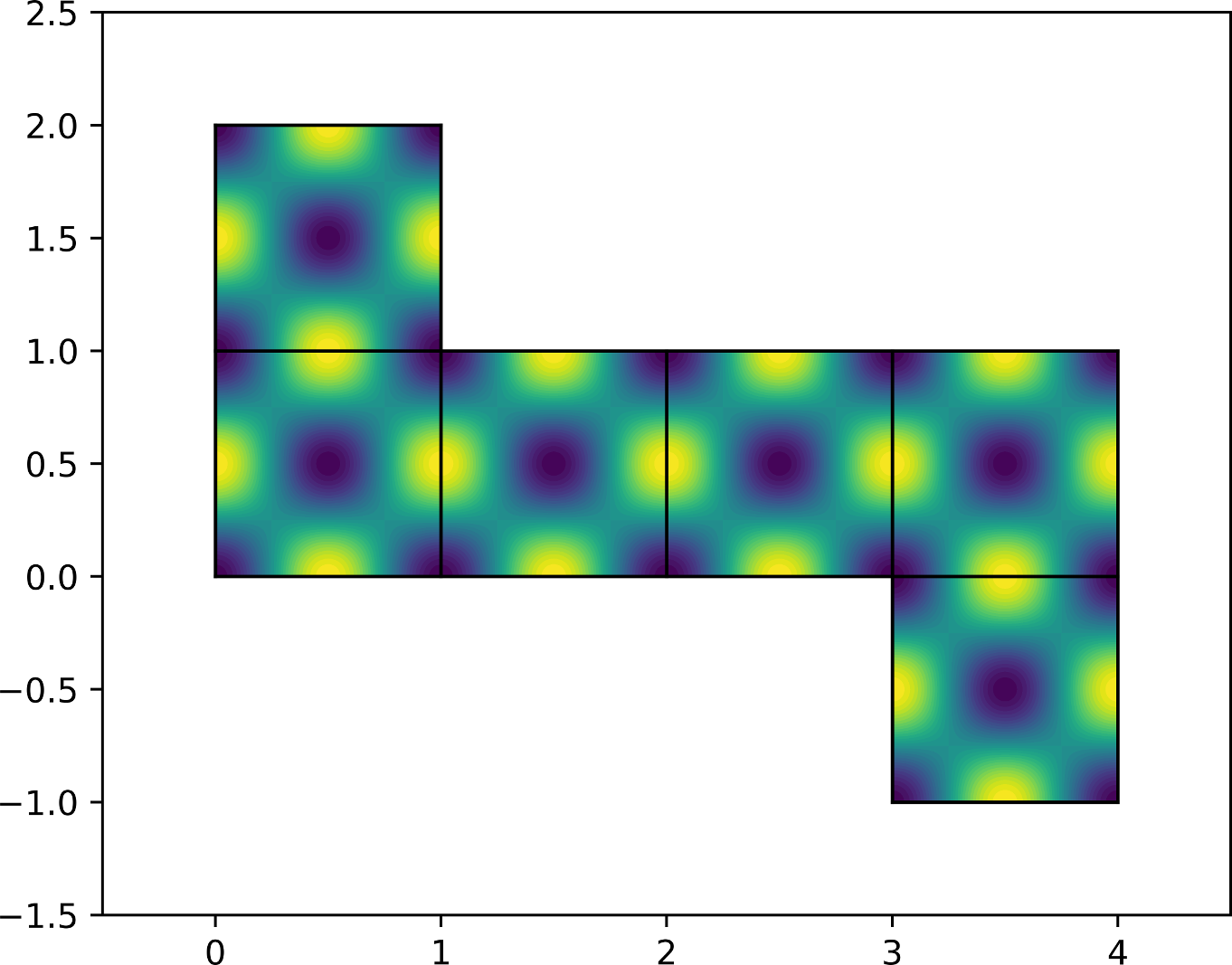}
    }
    \hfill
    \subfigure[For eigenvalue $\lambda_{38}=8.06$ (singular)]
    {
        \includegraphics[width=2.2in]{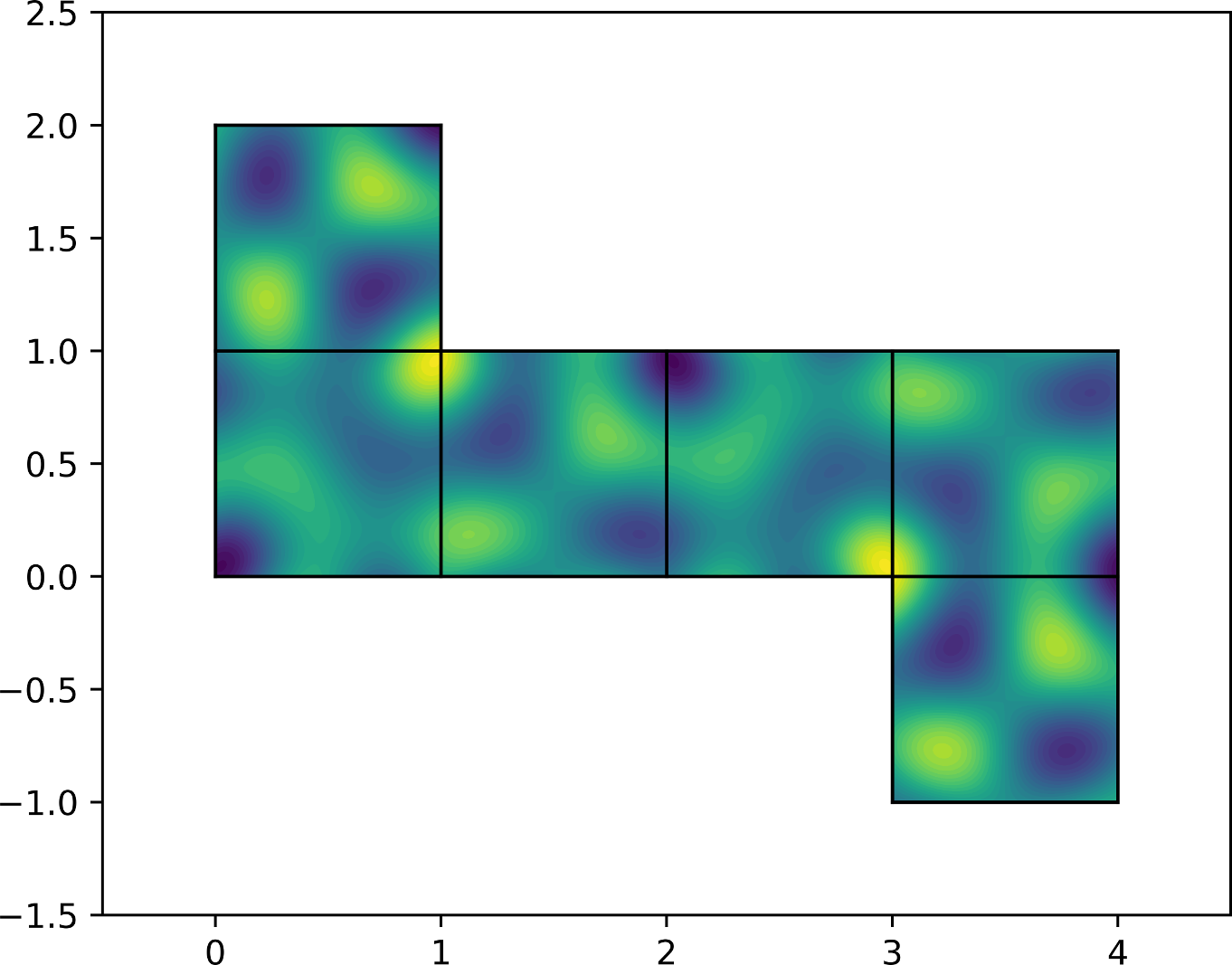}
    }
    \caption{Eigenfunctions on the cube}
    \label{fig:Figure 9}
\end{figure}

\vspace{-0.5cm}

\begin{figure}[h]
    \centering
    \subfigure[For eigenvalue $\lambda_{37}=8$ (nonsingular)]
    {
        \includegraphics[width=2.2in]{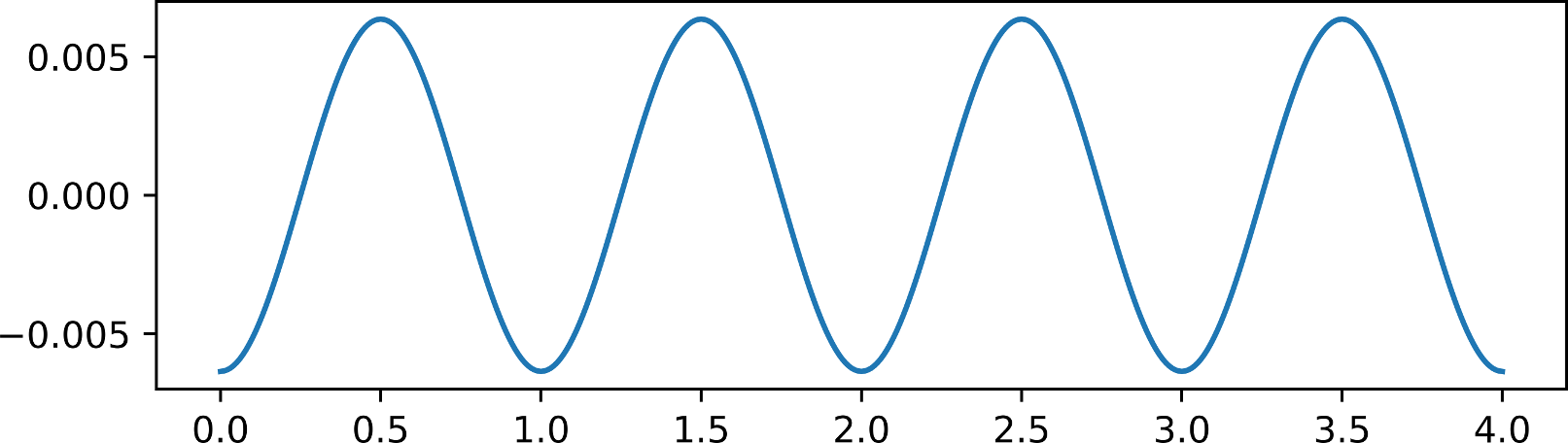}
    }
    \hfill
    \subfigure[For eigenvalue $\lambda_{38}=8.06$ (singular)]
    {
        \includegraphics[width=2.2in]{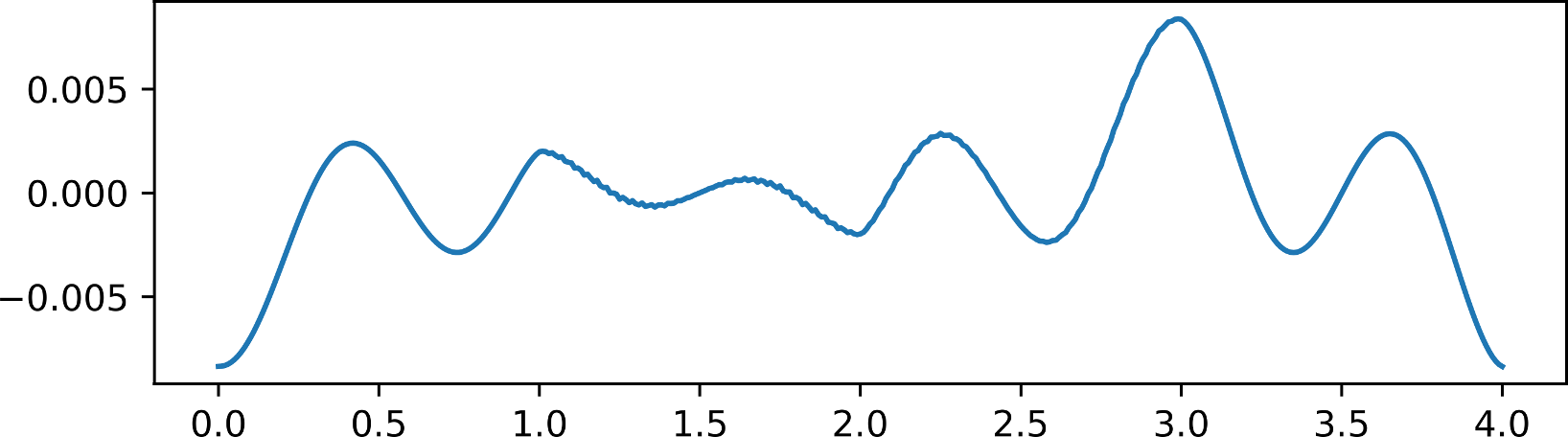}
    }
    \caption{Restriction of eigenfunctions on the cube to $y=0$}
    \label{fig:Figure 9}
\end{figure}

\vspace{-0.5cm}

\begin{figure}[h]
    \centering
    \subfigure[Eigenfunction from the 2-dim. eigenspace to eigenvalues $\lambda_{18}=\lambda_{19}=5.333$ (rotated)]
    {
        \includegraphics[width=2.2in]{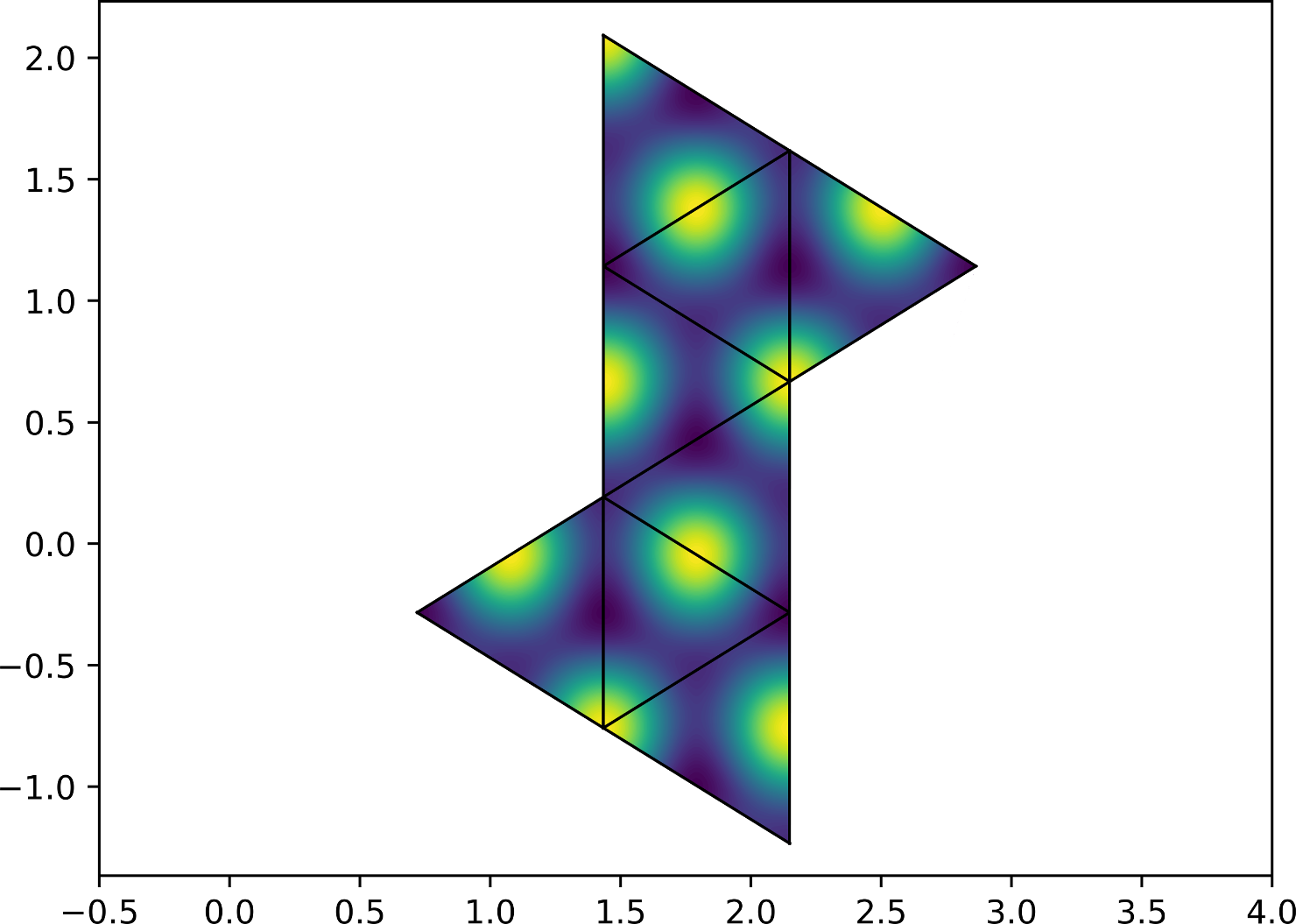}
    }
    \hfill
    \subfigure[Eigenfunction from the 2-dim. eigenspace to eigenvalues $\lambda_{56}=\lambda_{57}=16$]
    {
        \includegraphics[width=2.2in]{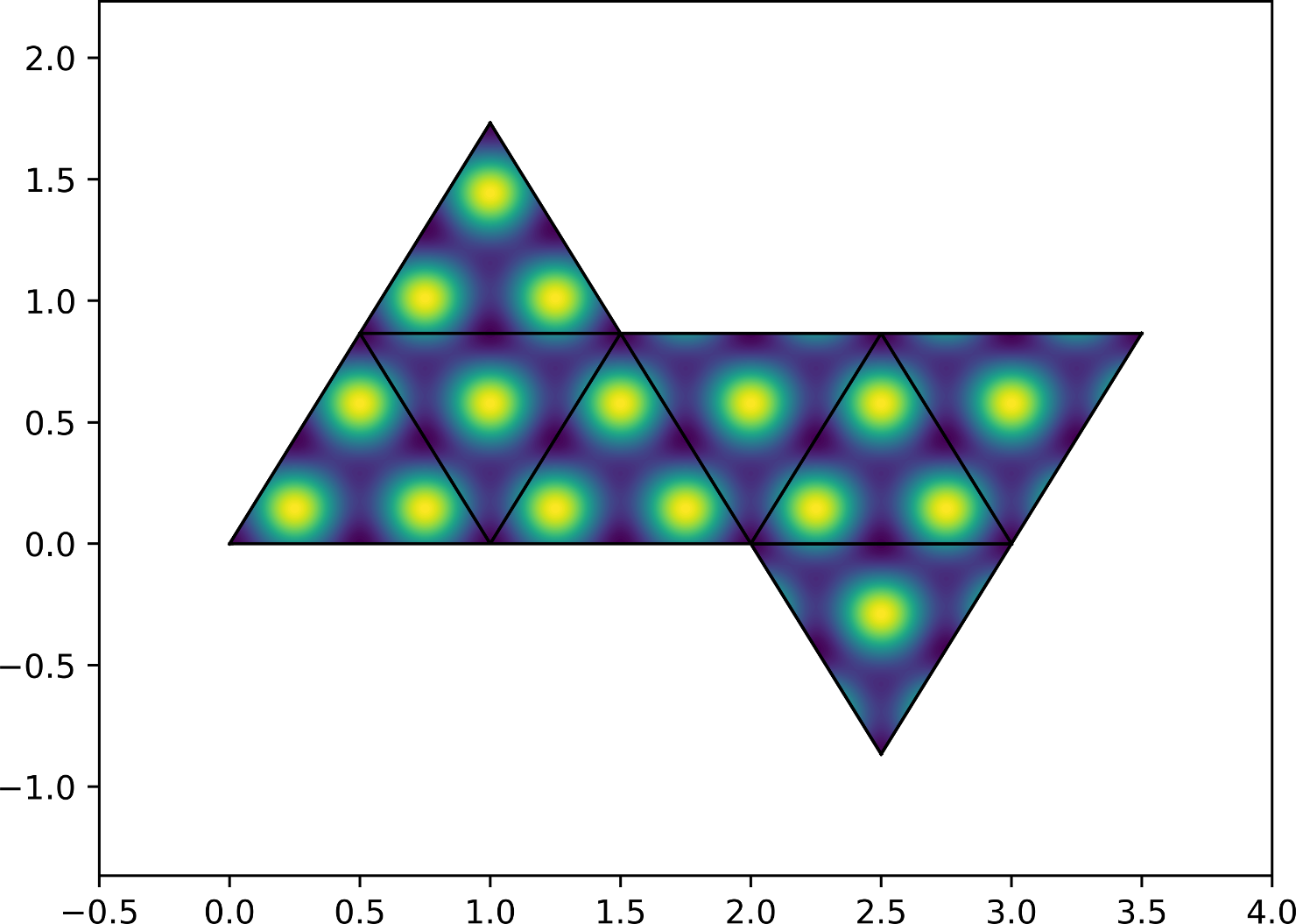}
    }
    \caption{Eigenfunctions on the octahedron:  (a) is obtained from (b) by rotating and dilating by $\sqrt{3}$, which reduces the eigenvalue by a factor of $\frac{1}{3}$}
    \label{fig:Figure 9}
\end{figure}

\begin{figure}[h]
\centering
\includegraphics[width=2.5in]{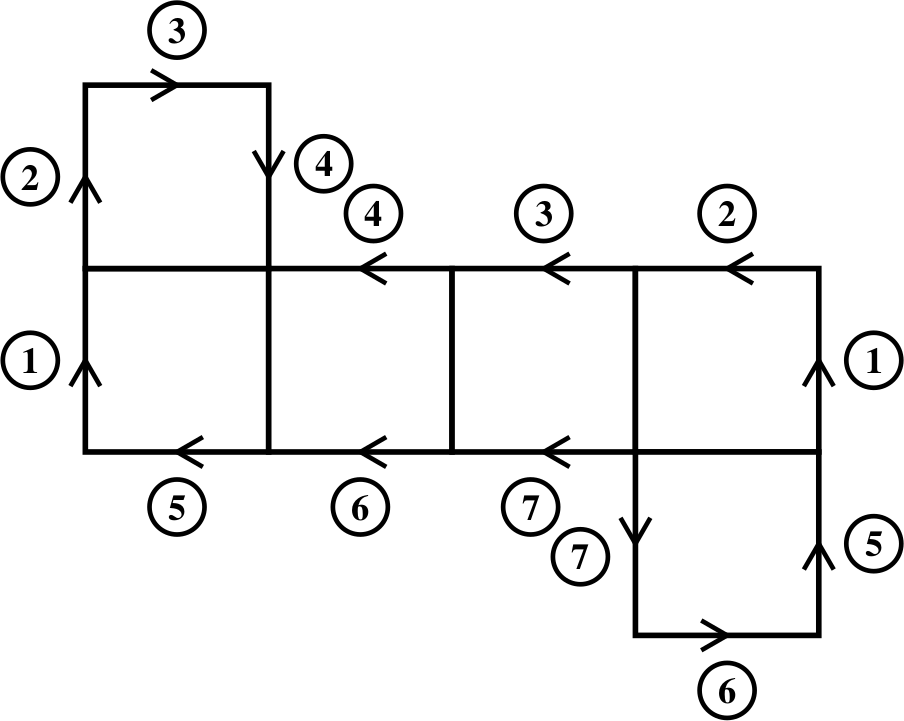}
\caption{Identified edges on the cube}
\label{}
\end{figure}

It is known that the eigenvalues form an increasing sequence $0 = \lambda_0 < \lambda_1 \leq \lambda_2 \leq \ldots$ tending to infinity, and satisfying the Weyl asymptotic law

\begin{equation}
\label{eq:weyl}
N(t) = \#\{\lambda_j \leq t\} \sim \frac{A}{4\pi}t
\end{equation}

\noindent where $A$ is the area of the surface. Also, the eigenfunctions may be sorted into symmetry types according to the irreducible representations of the symmetry group of the polyhedron. The eigenfunctions are real analytic functions in the complement of the singular vertex points. In particular, an eigenfunction is uniquely determined by its values on a single face. We will sort the eigenfunctions into two types: \textit{nonsingular} eigenfunctions have removable singularities at the vertices, while those with nonremovable singularities will be called \textit{singular}. We will see that  on the tetrahedron, all eigenfunctions are nonsingular, while for the other regular polyhedra both types exist, with most falling into the singular type. The nonsingular eigenfunctions may all be expressed in terms of trigonometric polynomials (except for the dodecahedron, because the eigenfunctions of the Laplacian on a regular pentagon do not have such representation), and we will give explicit formulas for them. The singular eigenfunctions constitute new families of ``special functions''.

We will begin our discussion in section $2$ with the tetrahedron. Here the story is very simple, because there exists a two-fold covering by a hexagonal torus (ignoring singular points). This allows us to lift all eigenfunctions on the tetrahedron to the torus, where everything is well known. One explanation for the simplicity in this case is that the cone angles are equal to $\pi$, which evenly divides $2\pi$. On the other hand, the symmetry group of the tetrahedron ($S_4$) has only two one-dimensional representations. A simple explanation for this is that the central reflection is not a symmetry of the tetrahedron. What this means is that the only possible symmetry behaviors for a one-dimensional eigenspace is either total symmetry or total skew-symmetry with respect to all reflections.

Our next example, discussed in section $3$, is the octahedron. This is perhaps the most interesting example. Note that the central reflection is a symmetry of the octahedron and generates a central $\mathbb{Z}_2$ subgroup of the full symmetry group, which is $S_4 \times \mathbb{Z}_2$. We note that there are two different types of reflection symmetries. There are what we call \textit{in-face} reflections, that reflect in the face bisectors of two opposite faces, permuting the other six faces, and the \textit{face-to-face} reflections that reflect in the boundary edges of faces and permute adjacent faces. This allows four different types of one-dimensional eigenspaces, with symmetry that we denote $\pm\pm$, the first $\pm$ indictating symmetry or skew-symmetry with respect to in-face reflections, and the second $\pm$ indicating symmetry or skew-symmetry with respect to face-to-face reflections. The $++$ and $--$ symmetry types coincide with the tetrahedron eigenspaces: every one-dimensional tetrahedron eigenfunction restricted to a face is also an octahedron eigenfunctions restricted to a face, and then reflecting evenly ($++$) or oddly ($--$) to extend to the whole surface. But there are also one-dimensional non-singular eigenspaces corresponding to the $+-$ and $-+$ symmetry types. Coincidentally they correspond to the same eigenvalue as one of the tetrahedron eigenfunctions. They are also given by trigonometric polynomials, but they are distinct from the tetrahedron eigenfunctions.

But that is not the end of the story for nonsingular eigenfunctions on the octahedron. For each of the above one-dimensional eigenspaces, there is an ``englargement'' by a factor of three in area. We take the eigenfunction restricted to one-sixth of a face, rotate and dilate it to obtain a function on a half-face, and then reflect in a pattern determined by the symmetry type. We then obtain an eigenfunction associated to the eigenvalue $\frac{\lambda}{3}$ where $\lambda$ is the eigenvalue for the original eigenfunction. These two types of constructions appear to exhaust the class of nonsingular eigenfunctions.

In sections $4$ and $5$ we discuss the icosahedron and the cube. We omit a discussion of the dodecahedron. Although it has nonsingular eigenfunctions, they are related to eigenfunctions on a regular pentagon, and so are not represented by trigonometric polynomials and do not have recognizable eigenvalues. We leave the study of these to the future, as well as the study of higher dimensional analogs.

In addition to the Weyl asymptotic law (1.2) for the eigenvalue counting function, we will examine a more refined estimate where we add a constant $c$ depending on the polyhedron to $\frac{A}{4\pi}t$. The difference

\begin{equation}
D(t) = N(t) - \left(\frac{A}{4\pi}t + c\right)
\end{equation}

only shows its advantages when averaged:

\begin{equation}
A(t) = \frac{1}{t}\int_0^t D(s)\text{ ds.}
\end{equation}

In the case of the tetrahedron, the results of [1] show that the choice $c=\frac{1}{2}$ implies

\begin{equation}
A(t)=O(t^{-\frac{1}{4}})
\end{equation}

and more precisely

\begin{equation}
g(t)=t^{\frac{1}{2}}A(t^2)
\end{equation}

is asymptotically equal to a uniformly almost periodic function of mean value 0. For the other three polyhedra, there are conjectures  in [5] for the correct choice of the constant that should lead to the same estimate. We support this conjectures by presenting experimental evidence, calculated using the Finite Element Method (FEM). We use linear splines and regular meshes with different refinements. For example, resolution 128 means, that each side of each face contains 128 vertices, and overall, that the mesh for the tetrahedron contains 33153 vertices, the octahedron 66177, the icosahedron 165249, and the cube 99201. We will extrapolate the data by taking the limit of an exponential fit using resolutions 32, 64 and 128, and compare the extrapolated data with the correct values for the tetrahedron.

The website \url{http://pi.math.cornell.edu/\~polyhedral} contains the programs we used and more data on eigenvalues and eigenfunctions.

For readers who are not familiar with the theory of irreducible representations of finite groups, it is possible to read this paper by skipping the passages that deal with this topic. Interested readers may find a nice introduction to the topic in [2] or [3]. There is a rather persuasive ideology that asserts that if a mathematical problem possesses symmetry, then it is highly worthwhile to try to exploit the symmetry in attacking the problem. We will be pleased if we can recruit some readers to embrace this point of view.

\newpage

\section{The Tetrahedron.}

The surface of a regular tetrahedron may be represented as a planar region with boundary identifications, as shown in Figure 5.

\begin{figure}[h]
\centering
\includegraphics[width=2.0in]{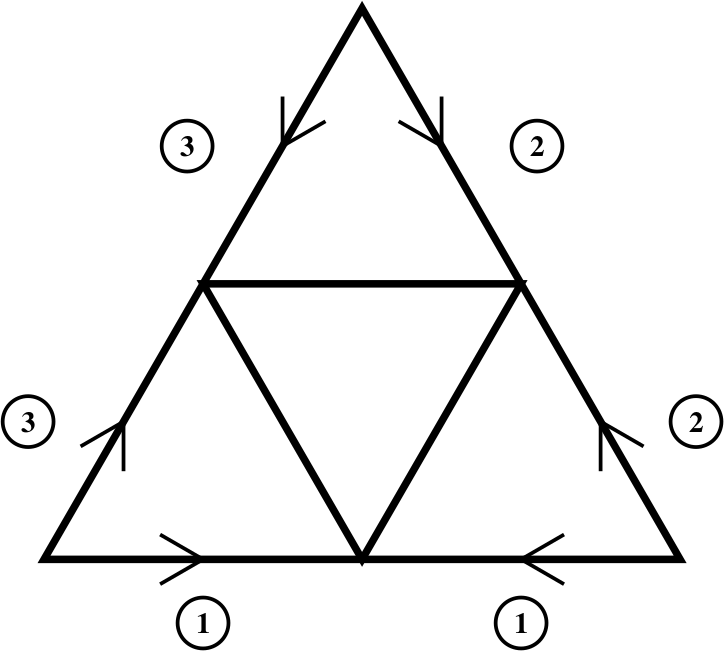}
\caption{Identified edges on the tetrahedron}
\label{fig:Figure 2}
\end{figure}

All eigenfunctions on the tetrahedron are nonsingular. We show two examples in Figure $6$ and the respective graphs of the restriction to line segments passing through vertices in Figure $7$.

\begin{figure}[h]
    \centering
    \subfigure[For eigenvalue $\lambda_{121}=64$ (nonsingular)]
    {
        \includegraphics[width=2.2in]{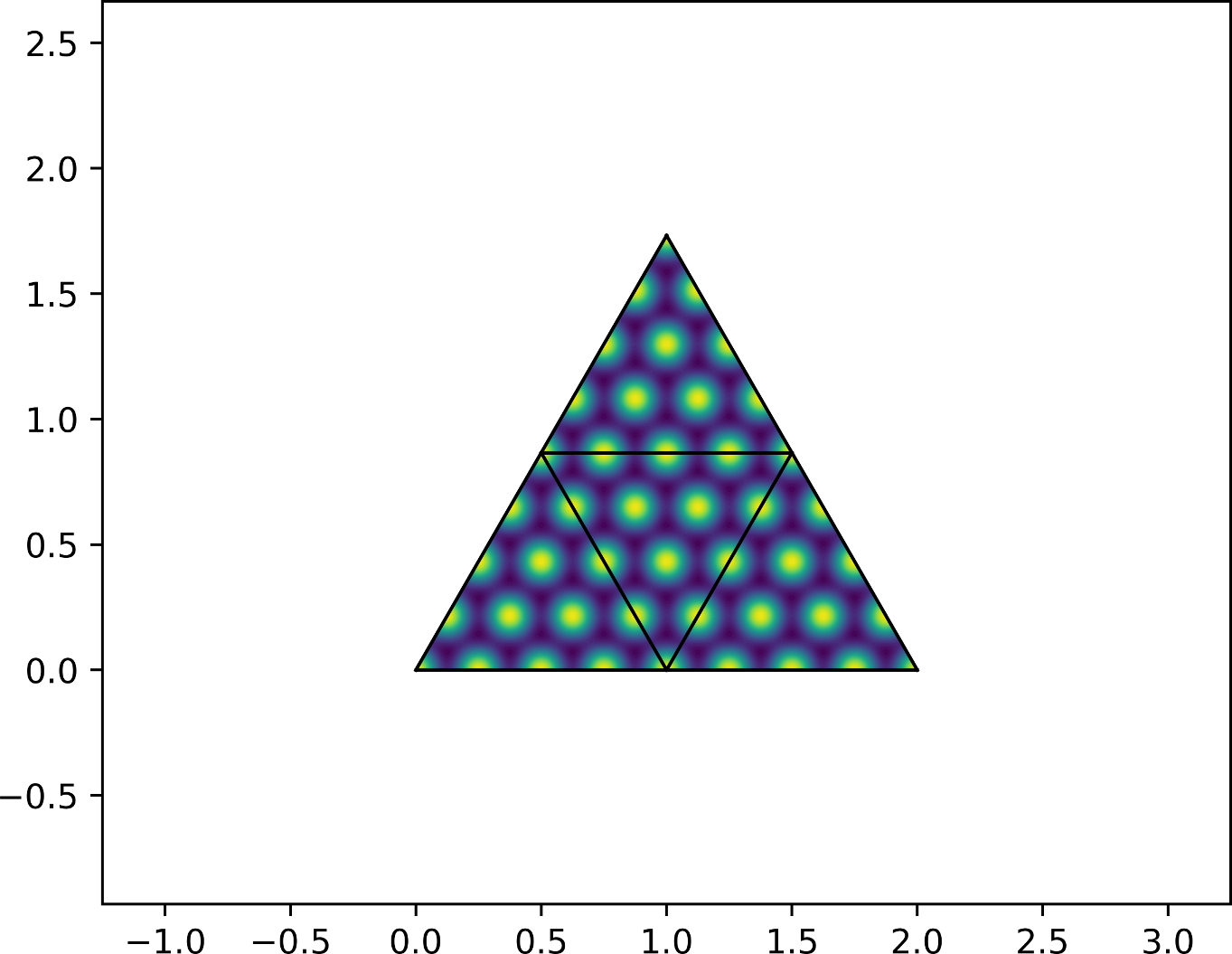}
        \label{fig:Figure1a}
    }
    \subfigure[For eigenvalue $\lambda_{122}=67$ (singular)]
    {
        \includegraphics[width=2.2in]{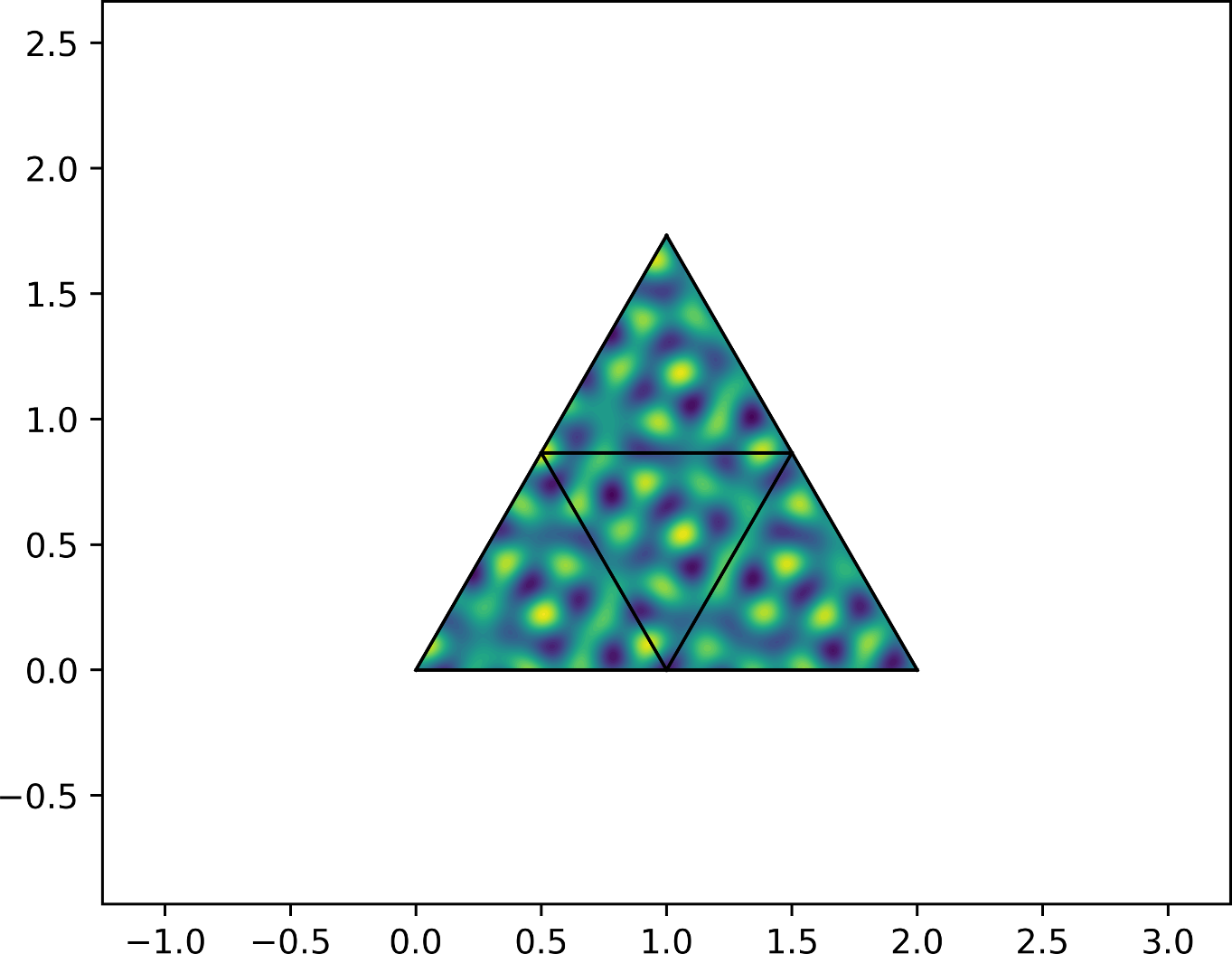}
        \label{fig:Figure2b}
    }
    \caption{Eigenfunctions on the tetrahedron}
    \label{fig:Figure 9}
	\vspace{0.5cm}
    \subfigure[For eigenvalue $\lambda_{121}=64$ (nonsingular)]
    {
        \includegraphics[width=2.2in]{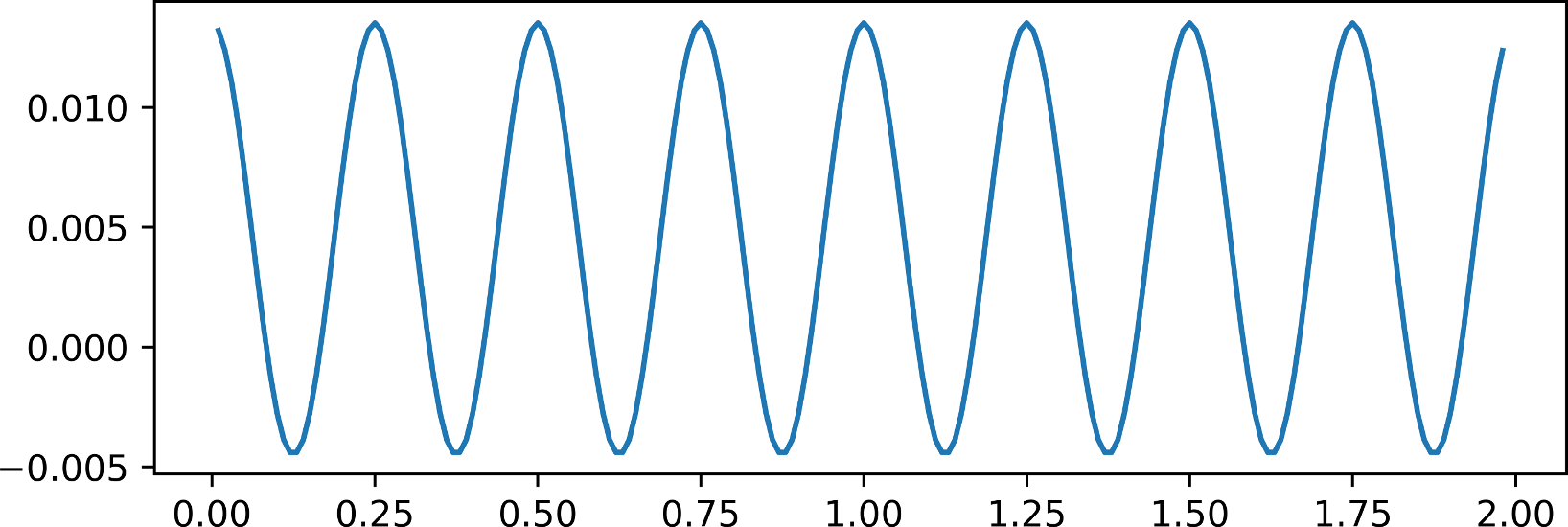}
        \label{fig:Figure1a}
    }
    \subfigure[For eigenvalue $\lambda_{122}=67$ (singular)]
    {
        \includegraphics[width=2.2in]{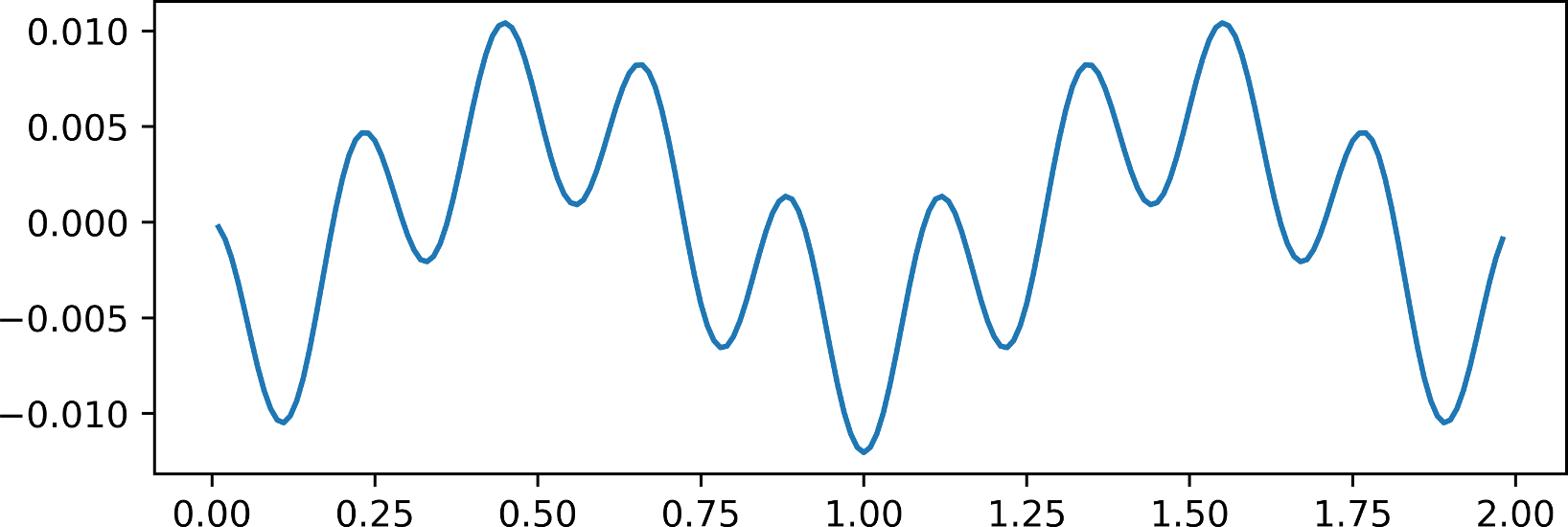}
        \label{fig:Figure2b}
    }
    \caption{Restriction of eigenfunctions on the tetrahedron to $y=0$}
    \label{fig:Figure 9}
\end{figure}

\newpage

A simple observation is that if we delete the vertices, there is a two-fold covering by a torus, again with singular points deleted, as shown in Figure $8$.

\begin{figure}[h]
\centering
\includegraphics[width=2in]{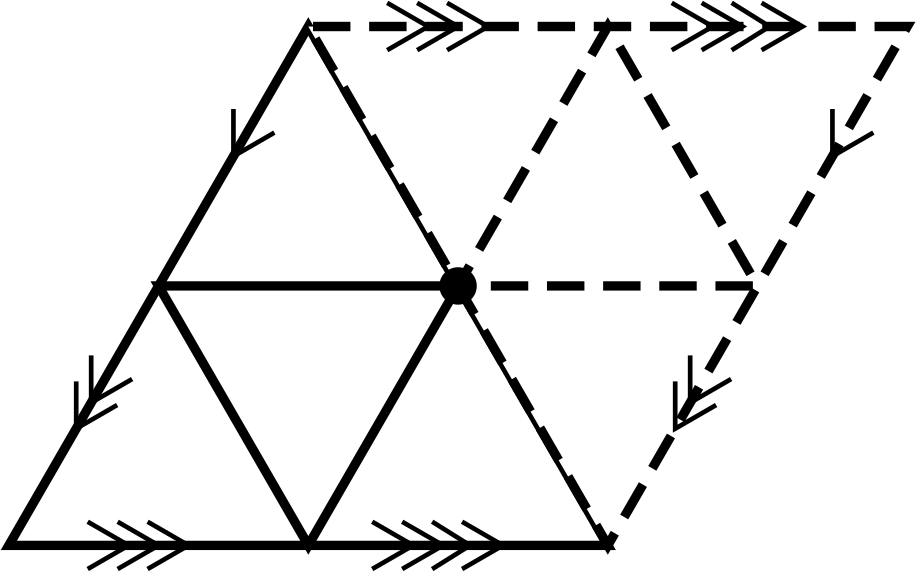}
\caption{The torus covering}
\label{fig:Figure 3}
\end{figure}

The covering identifications are obtained by a half ($180^\circ$) rotation $R$ about the indicated center point. In other words, any smooth function on the tetrahedron lifts to a smooth function on the torus that is invariant under $R$, 

\begin{equation}
u(Rx) = u(x), 
\end{equation}

\noindent and, conversely, every smooth function on the torus satisfying (2.1) defines a smooth function on the tetrahedron. In particular, this applies to eigenfunctions of the Laplacian. Since points are removable singularities for the Laplacian in the plane, the eigenfunctions on the torus with singular points removed are actually eigenfunctions on the torus, and these are all known.

To be specific, we will assume that the edges on the tetrahedron all have length one. Then a single face may be chosen to have vertices $(0,0),(1,0), (\frac{1}{2},\frac{\sqrt{3}}{2})$. The torus is the span of $a(2,0) + b(1,\sqrt{3})$ for $0 \leq a \leq 1$, $0 \leq b \leq 1$. This means that a function on the torus may be lifted to a periodic function on the plane with periodicity conditions

\begin{align*}
u\left(x+(2,0)\right) &= u(x) \\
u(x+(1,\sqrt{3})) &= u(x)
\end{align*}

The dual lattice to the lattice generated by $(2,0)$ and $(1,\sqrt{3})$ has generators 
\begin{equation*}
\overset{\rightarrow}{u} = \left(\frac{1}{2}, \frac{\sqrt{3}}{6}\right)\text{ and }\overset{\rightarrow}{v} = \left(0,\frac{\sqrt{3}}{3}\right)\text{,}
\end{equation*}

\noindent so all eigenfunctions on the torus are generated by 

\begin{equation*}
e^{2\pi ix \cdot (k\overset{\rightarrow}{u} + j\overset{\rightarrow}{v})}
\end{equation*}

\vspace{0.2cm}

\noindent for $(k,j) \in \mathbb{Z}^2$, with eigenvalue $\frac{4\pi^2}{3}\left(j^2+k^2+jk\right)$. The condition (2.1) then gives the tetrahedron eigenfunctions as 

\begin{equation}
\label{eq:torfuncs}
\cos\left(2\pi x \cdot (k\overset{\rightarrow}{u} + j\overset{\rightarrow}{v})\right) = \cos\left(2\pi\left(\frac{1}{2}kx_1+\left(\frac{\sqrt{3}}{3}j+\frac{\sqrt{3}}{6}k\right)x_2\right)\right)\text{.}
\end{equation}

\vspace{0.1cm}

The symmetry group of the tetrahedron is isomorphic to the permutation group $S_4$ on $4$ letters, since the $4$ vertices of the tetrahedron may be arbitrarily permuted by an isometry. The group is generated by reflections in a bisector of a face, that will also reflect in a bisector of a second face and switch the other two faces by reflection in their common edge (Figure $9$).

\begin{figure}[h]
\centering
\includegraphics[width=1.8in]{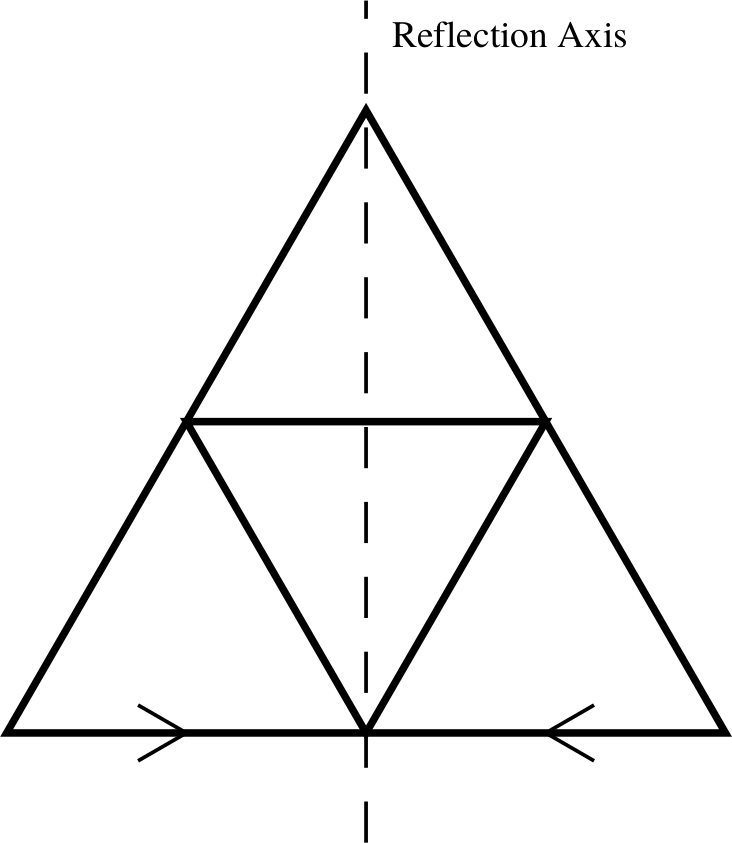}
\caption{Reflection symmetry}
\label{fig:Figure 4}
\end{figure}

This group has five distinct irreducible representations of dimensions $3,3,2,1,1$. In particular, functions that transform according to one of the 1-dimensional representations (denoted by $1^+$ and $1^-$) are easily described as being symmetric $(1^+)$ and skew-symmetric $(1^-)$ respectively, with respect to all reflections. Such functions are uniquely determined by their restriction to a fundamental domain, one-sixth of a face, and extended accordingly (see Figure $10$ for the skew-symmetric case).

\begin{figure}[h]
\centering
\includegraphics[width=2in]{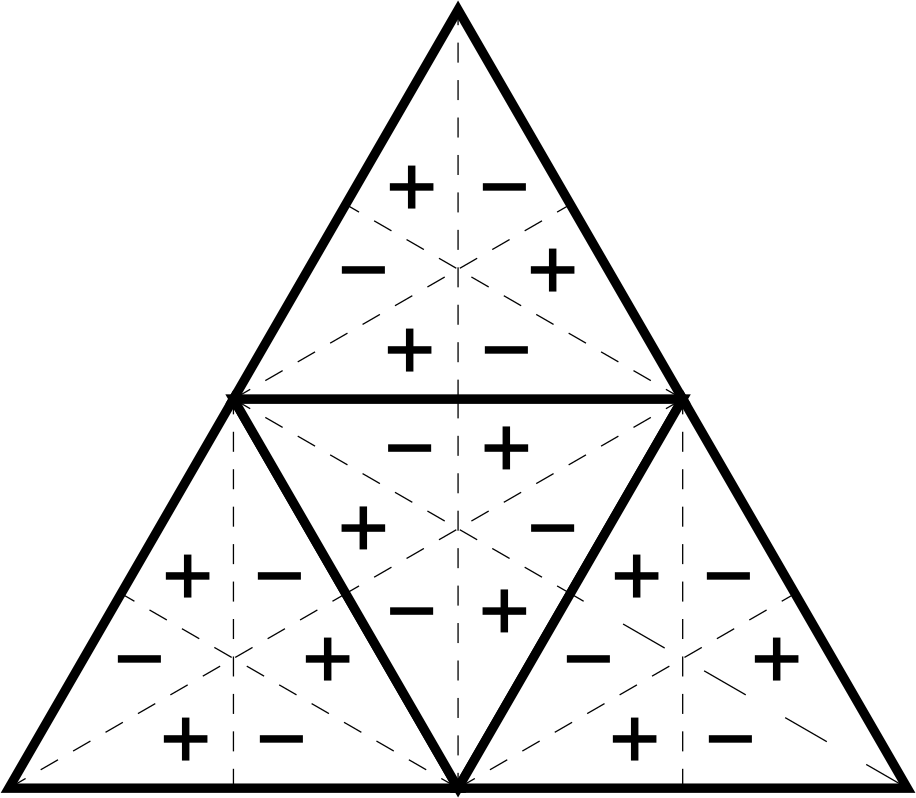}
\caption{Skew-symmetric (1-) reflection}
\label{fig:Figure 5}
\end{figure}

It is clear from this description that $1^+$ and $1^-$ type eigenfunctions are actually well defined on the smaller torus consisting of two faces with edges identified, as the torus translations are expressible as the product of two reflections. In terms of the respresentation as sums of functions of the form (2.2) this just means that the integers $k$ and $j$ must both be even. To describe all eigenfunctions with these symmetries we look at the points in the lattice generated by $\overset{\rightarrow}{u}$ and $\overset{\rightarrow}{v}$, the hexagonal lattice, and group them in concentric hexagons around the origin, Figure $11$.

\begin{figure}[h]
\centering
\includegraphics[width=1.6in]{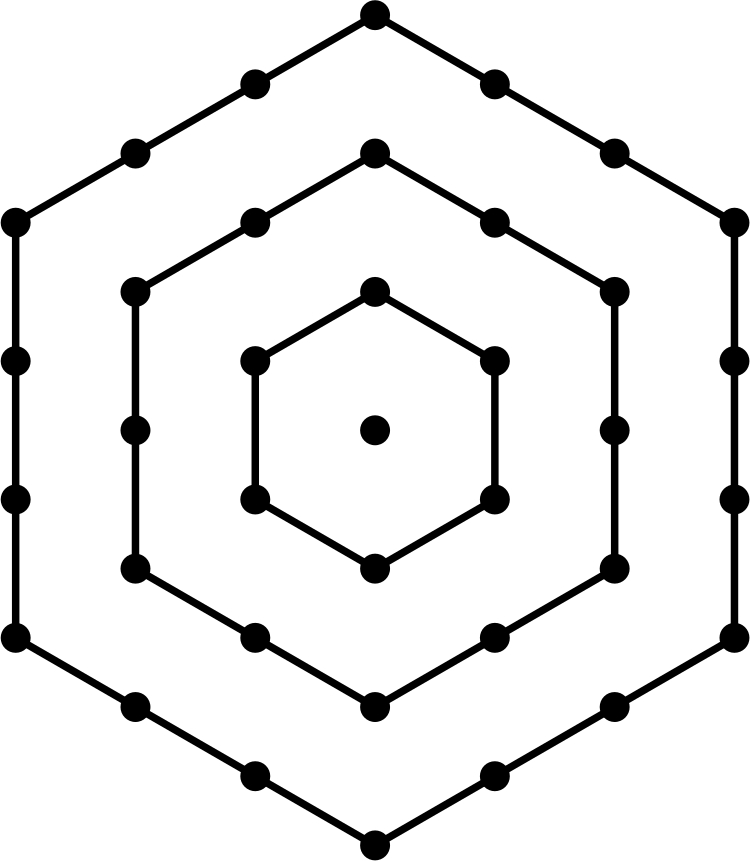}
\caption{The hexagonal lattice}
\label{fig:Figure 6}
\end{figure}

The dihedral-$6$ group of symmetries of the hexagon is generated by $6$ reflections, and the orbits under this action have generically $12$ elements. There are also two types of $6$-element orbits: type $(1)$ consists of the vertices of a concentric hexagon, and type $(2)$ consists of the midpoints. Note that type $(2)$ occurs only in every other concentric hexagon (the origin is a singular orbit). The $1^+$ and $1^-$ eigenfunctions correspond to sums of (2.2) over an orbit, with alternating $+$ and $-$ signs for the $1^-$ case. The $6$-element orbits do not have $1^-$ eigenfunctions. The distribution of signs of the $1^-$ eigenfunction on a generic orbit is shown in Figure $12$, which also gives the $(k,j)$ values at each point. 

\vspace{0.2cm}

\begin{figure}[h]
\centering
\includegraphics[width=2.7in]{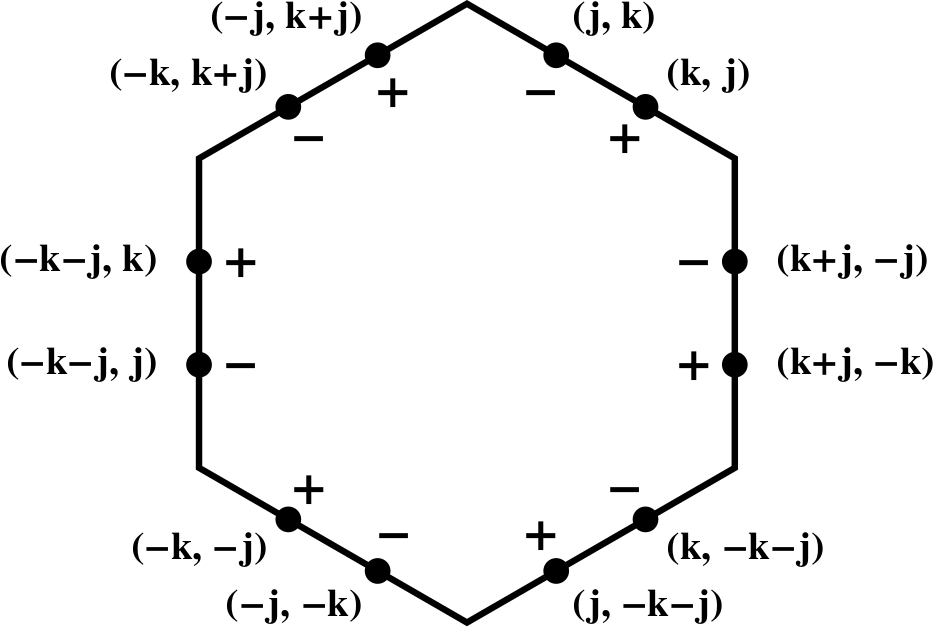}
\caption{A generic orbit ($j > k > 0$).}
\label{fig:Figure 7}
\end{figure}

\noindent Note that diametrically opposed vertices have the same sign, which is consistent with the representation as cosines. From this we can read off the explicit form of the eigenfunctions $u_+$ and $u_-$ associated to a generic orbit:

\begin{equation}
\begin{aligned}
\label{eq:uplusmin}
u_{\pm } &= \cos\left(2\pi x \cdot \left(k\overset{\rightarrow}{u} + j\overset{\rightarrow}{v}\right)\right) + \cos\left(2\pi x \cdot \left(\left(k+j\right)\overset{\rightarrow}{u} - k \overset{\rightarrow}{v}\right)\right)
+ \cos\left(2\pi x \cdot \left(j\overset{\rightarrow}{u} - \left(k+j\right)\overset{\rightarrow}{v}\right)\right) \\
 &\pm \left(\cos\left(2\pi x \cdot \left(j\overset{\rightarrow}{u}+k\overset{\rightarrow}{v}\right)\right) + \cos\left(2\pi x \cdot \left(\left(k+j\right)\overset{\rightarrow}{u} - j \overset{\rightarrow}{v}\right)\right)
 + \cos\left(2\pi x \cdot \left(k\overset{\rightarrow}{u} - \left(k+j\right)\overset{\rightarrow}{v}\right)\right)\right)
\end{aligned}
\end{equation}

\vspace{0.1cm}

with $k > j > 0$ and both even.

The type $1$ nongeneric orbit corresponds to $j=0$. In that case the six terms collapse to the first three:

\begin{equation}
\label{eq:uplus}
u_{+} = \cos(2\pi x \cdot \overset{\rightarrow}{u}) + \cos(2\pi k(\overset{\rightarrow}{u} - \overset{\rightarrow}{v})) + \cos(2\pi k x \cdot \overset{\rightarrow}{v})
\end{equation}

\vspace{0.1cm}

The type $2$ nongeneric orbit corresponds to $j=k$. In that case the sum collapses to 

\begin{equation}
\label{eq:uplus2}
u_+ = \cos(2\pi k x \cdot (\overset{\rightarrow}{u} + \overset{\rightarrow}{v})) + \cos(2\pi k x \cdot (2\overset{\rightarrow}{u}-\overset{\rightarrow}{v})) + \cos(2\pi k x \cdot (\overset{\rightarrow}{u} - 2 \overset{\rightarrow}{v}))
\end{equation}

\vspace{0.1cm}

Of course the origin orbit corresponds to the constant eigenfunction. The general eigenvalues are all of the form $\frac{4}{3}\pi^2 N$ where $N$ is a nonnegative integer expressible as $j^2 + k^2 + jk$, with multiplicity equal to the number of distinct ways to express $N$ in this form with $k \geq 0$. Aside from the trivial zero eigenvalue, the multiplicities will typically be $3$ or $6$ (but occasionally higher due to coincidences). For eigenspaces containing a $1^+$ eigenfunction, or both $1^+$ and $1^-$ eigenspaces, the remaining orthogonal complement must transform according to the $2$-dimensional representation (one or two copies). For all the other eigenspaces ($k$ and $j$ not both even), those of multiplicity three must transform as one of the $3$-dimensional representations. We speculate that those of dimension six will split into a sum of each of the $3$-dimensional representations.

What is the proportion of $1^+$ and $1^-$ eigenfunctions among all eigenfunctions? A rough count will just include generic orbits. One quarter of them will have both $k$ and $j$ even, and among these orbits one sixth of each will be $1^+$ and $1^-$. Thus asymptotically $\frac{1}{24}$ of all eigenfunctions are of type $1^+$ and $1/24$ are of type $1^-$. Table 1 shows the first 60 eigenvalues normalized by dividing by $\frac{4}{3}\pi^2$. Since we know that these values are all integers, we can immidiately see the accuracy of our computations and of the extrapolation (Figure 13). We will continue using the corrected (integer) data.

\begin{figure}[h]
    \centering
    \subfigure[Eigenvalue data computed with different meshes, resolution 32 (blue), 64 (orange), 128 (green) and the corresponding extrapolation (red)]
    {
        \includegraphics[width=2.2in]{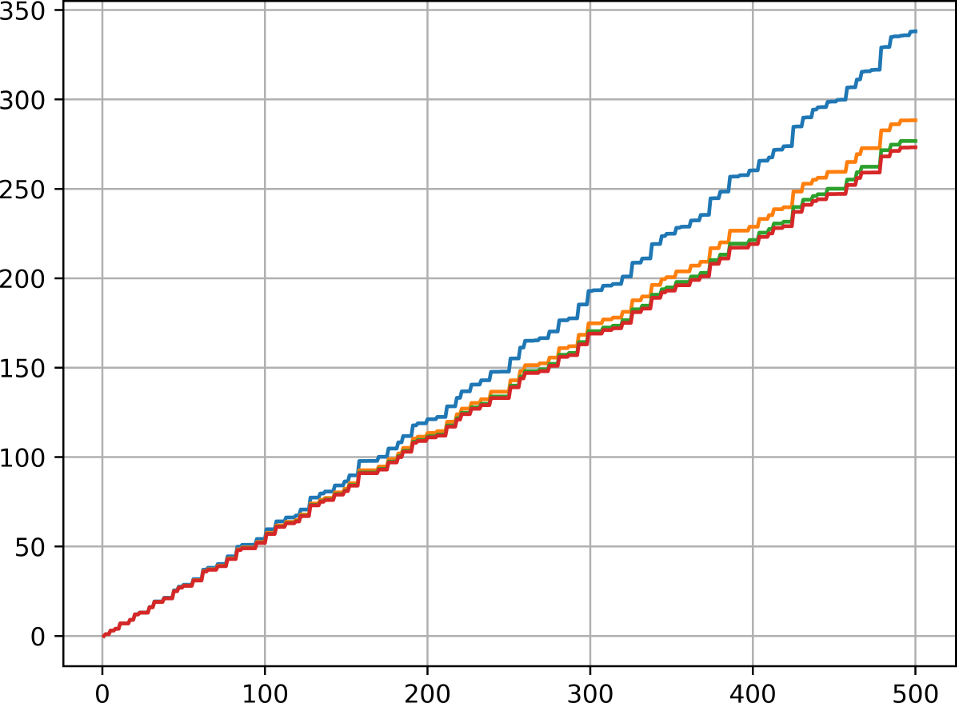}
        \label{fig:Figure1a}
    }
    \hfill
    \subfigure[Difference between res. 128 and the corrected data (blue), difference between the extrapolation and the corrected data (orange): the error between the corrected data and the computed data using the most refined mesh grows fast, but the extrapolated data is very accurate]
    {
        \includegraphics[width=2.2in]{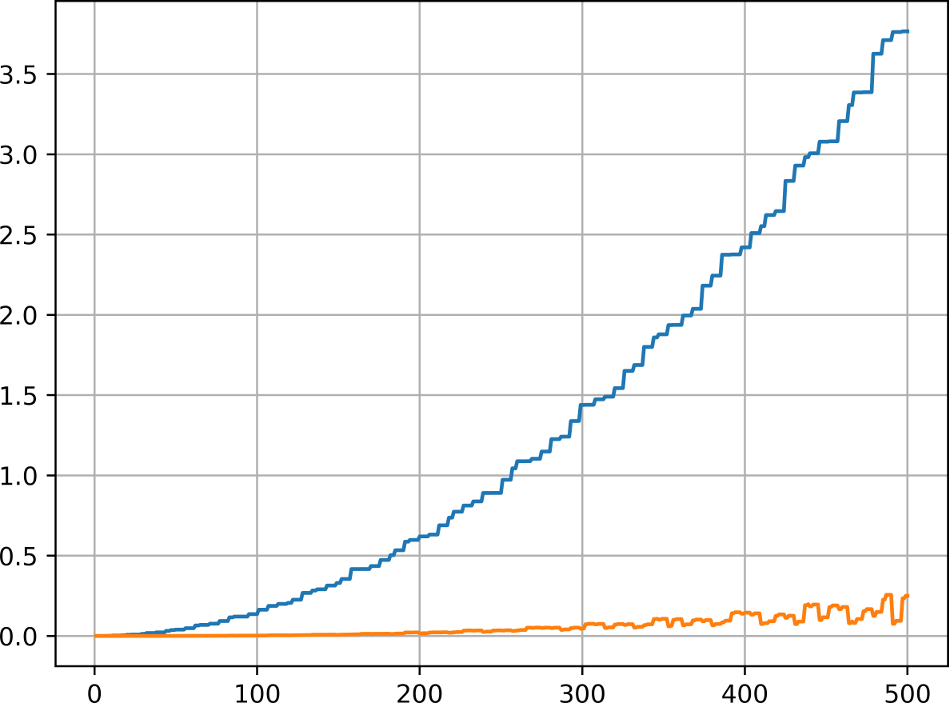}
        \label{fig:Figure2b}
    }
    \caption{Eigenvalues on the tetrahedron}
    \label{fig:Figure 9}
\end{figure}

\begin{table}[h!]
\centering
\makebox[0pt][c]{\parbox{\textwidth}{
\begin{minipage}[b]{0.32\hsize}
\centering
\begin{tabular}{| c | c |}
\hline
\# & Eigenvalue \\ 	
\hline 
  1 & 0\\
  2 & 1.00001 \\ 
  3 & 1.00001 \\
  4 & 1.00001 \\
  5 & 3.00011 \\
  6 & 3.00011 \\
  7 & 3.00011 \\
  8 & 4.00020 \\
  9 & 4.00020 \\
  10 & 4.00020 \\
  11 & 7.00062 \\
  12 & 7.00062 \\
  13 & 7.00062 \\
  14 & 7.00062 \\
  15 & 7.00062 \\
  16 & 7.00062 \\
  17 & 9.00102 \\
  18 & 9.00102 \\
  19 & 9.00102 \\
  20 & 12.00181 \\
\hline 
\end{tabular}
\end{minipage}
\hfill
\begin{minipage}[b]{0.32\hsize}
\centering
\begin{tabular}{| c | c |}
\hline
\# & Eigenvalue \\ 	
\hline 
  21 & 12.00181 \\
  22 & 12.00181 \\ 
  23 & 13.00212 \\
  24 & 13.00212 \\
  25 & 13.00212 \\
  26 & 13.00212 \\
  27 & 13.00212 \\
  28 & 13.00212 \\
  29 & 16.00321 \\
  30 & 16.00321 \\
  31 & 16.00321 \\
  32 & 19.00453 \\
  33 & 19.00453 \\
  34 & 19.00453 \\
  35 & 19.00453 \\
  36 & 19.00453 \\
  37 & 19.00453 \\
  38 & 21.00554 \\
  39 & 21.00554 \\
  40 & 21.00554 \\  
\hline 
\end{tabular}
\end{minipage}
\hfill
\begin{minipage}[b]{0.32\hsize}
\centering
\begin{tabular}{| c | c |}
\hline
\# & Eigenvalue \\ 	
\hline 
  41 & 21.00554\\
  42 & 21.00554 \\ 
  43 & 21.00554 \\
  44 & 25.00784 \\
  45 & 25.00784 \\
  46 & 25.00784 \\
  47 & 27.00915 \\
  48 & 27.00915 \\
  49 & 27.00915 \\
  50 & 28.00984 \\
  51 & 28.00984 \\
  52 & 28.00984 \\
  53 & 28.00984 \\
  54 & 28.00984 \\
  55 & 28.00984 \\
  56 & 31.01206 \\
  57 & 31.01206 \\
  58 & 31.01206 \\
  59 & 31.01206 \\
  60 & 31.01206 \\
\hline 
\end{tabular}
\end{minipage}
}}
\vspace{0.5cm}
\caption{Normalized eigenvalues on the tetrahedron, res. 128: we can see the error (the deviation from the integer value) growing}
\end{table}

\vspace{-0.4cm}

The eigenvalue counting function

\begin{equation}
N(t) = \#\{\lambda_j\leq t\}
\end{equation}

where $\lambda_j$ are the eigenvalues (counted according to multiplicity) is easy to understand in terms of the counting function $N_T$ for the covering torus, which is described explicitly in [1]. From Figure $11$ it is apparent that, aside from the constant eigenfunction, the eigenfunctions on $T$ split evenly among these that are symmetric and skew-symmetric with respect to R. So $N(t) = \frac{1}{2} N_T(t) + \frac{1}{2}$.

\newpage

In Figure $14$ we show the graphs of the eigenvalue counting function N(t), the difference

\begin{equation}
D(t) = N(t) - \left(\frac{\sqrt{3}}{4\pi}t + \frac{1}{2}\right)\text{,}
\end{equation}

the average of the difference

\begin{equation}
A(t) = \frac{1}{t}\int_0^t D(s)\text{ ds}
\end{equation}

\noindent and the average after rescaling

\begin{equation}
g(t)=t^{\frac{1}{2}}A(t^2)
\end{equation}

\noindent based on our numerical computations. It is shown in [1] Thm. 2.3 that $g(t)$ is asymptotic to a uniformly almost periodic function of mean value zero with an explicitly given Fourier series. The constant $\frac{\sqrt{3}}{4\pi}$ comes from the Weyl term.

\begin{figure}[h]
    \centering
    \subfigure[The counting function $N(t)$ (blue) and $\left(\frac{\sqrt{3}}{4\pi}t + \frac{1}{2}\right)$ (orange)]
    {
        \includegraphics[width=2.1in]{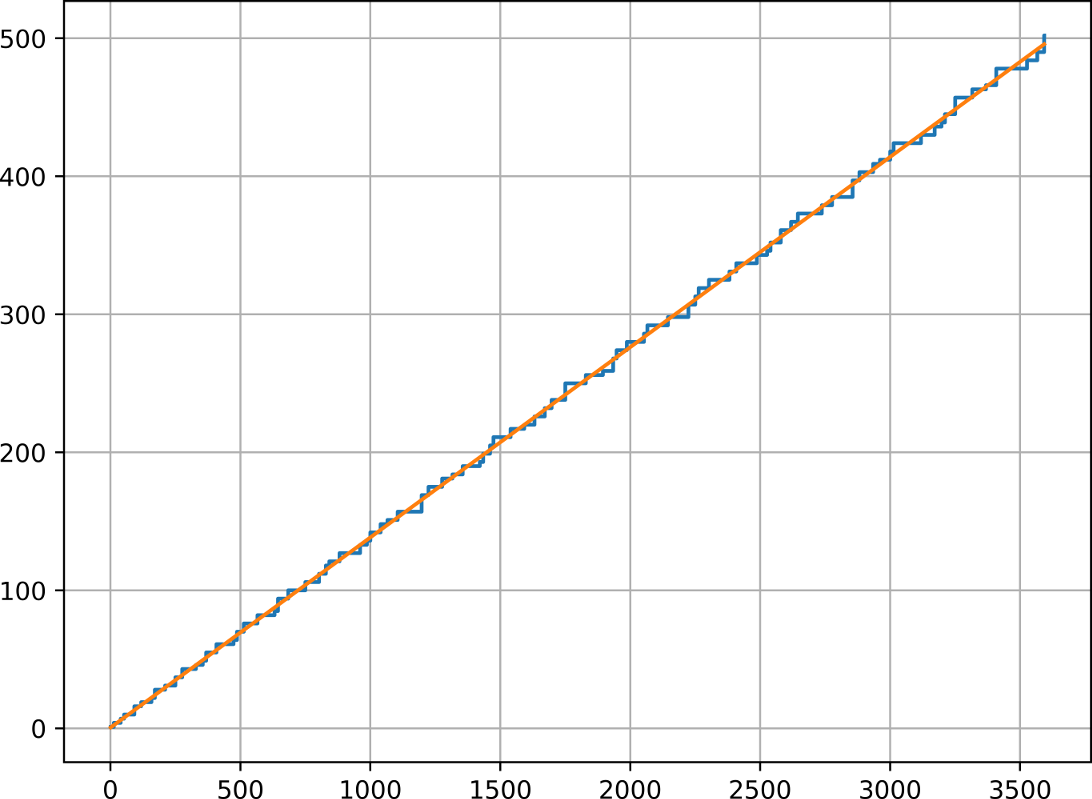}
        \label{fig:Figure1a}
    }
    \hfill
    \subfigure[The difference $D(t)$]
    {
        \includegraphics[width=2.1in]{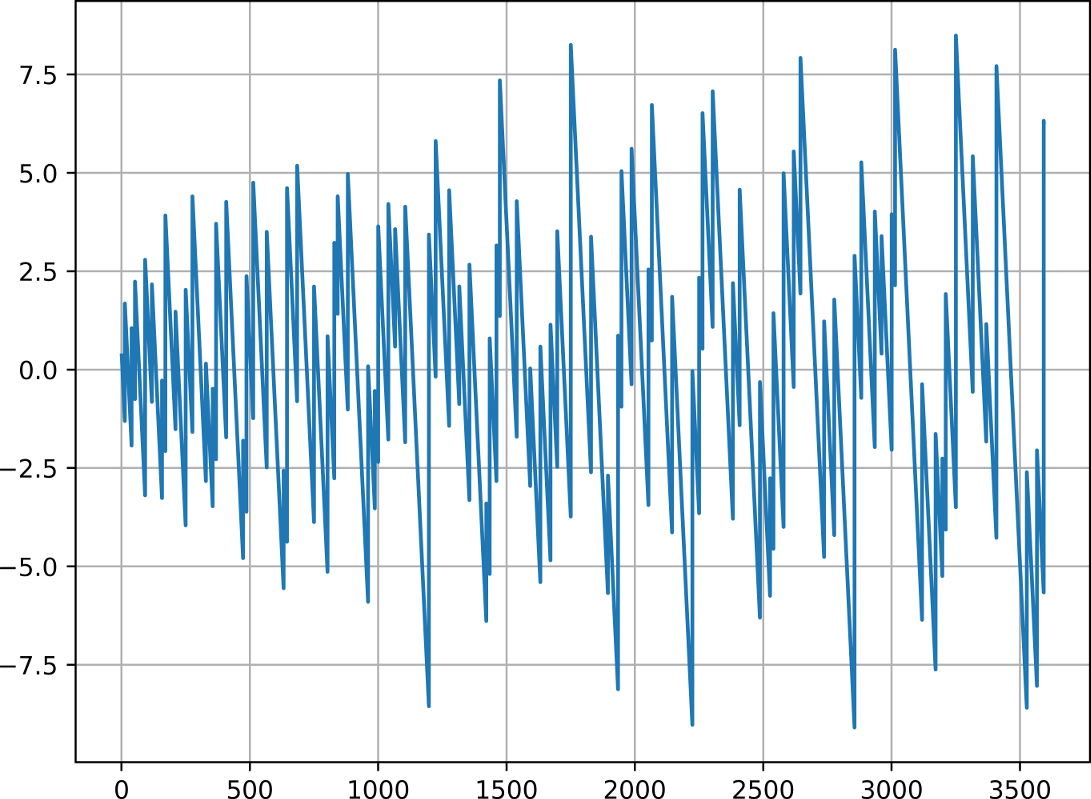}
        \label{fig:Figure2b}
    }
    \subfigure[The average of the difference $A(t)$]
    {
        \includegraphics[width=2.1in]{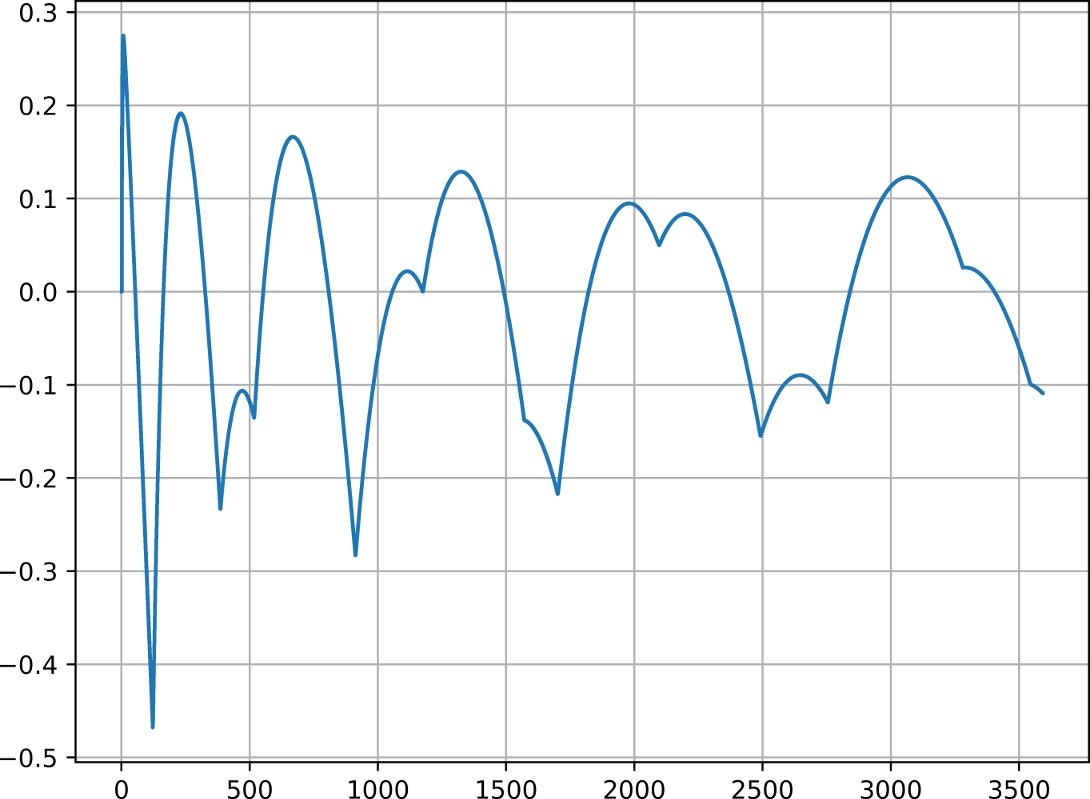}
        \label{fig:Figure1a}
    }
    \hfill
    \subfigure[The average after rescaling $g(t)$]
    {
        \includegraphics[width=2.1in]{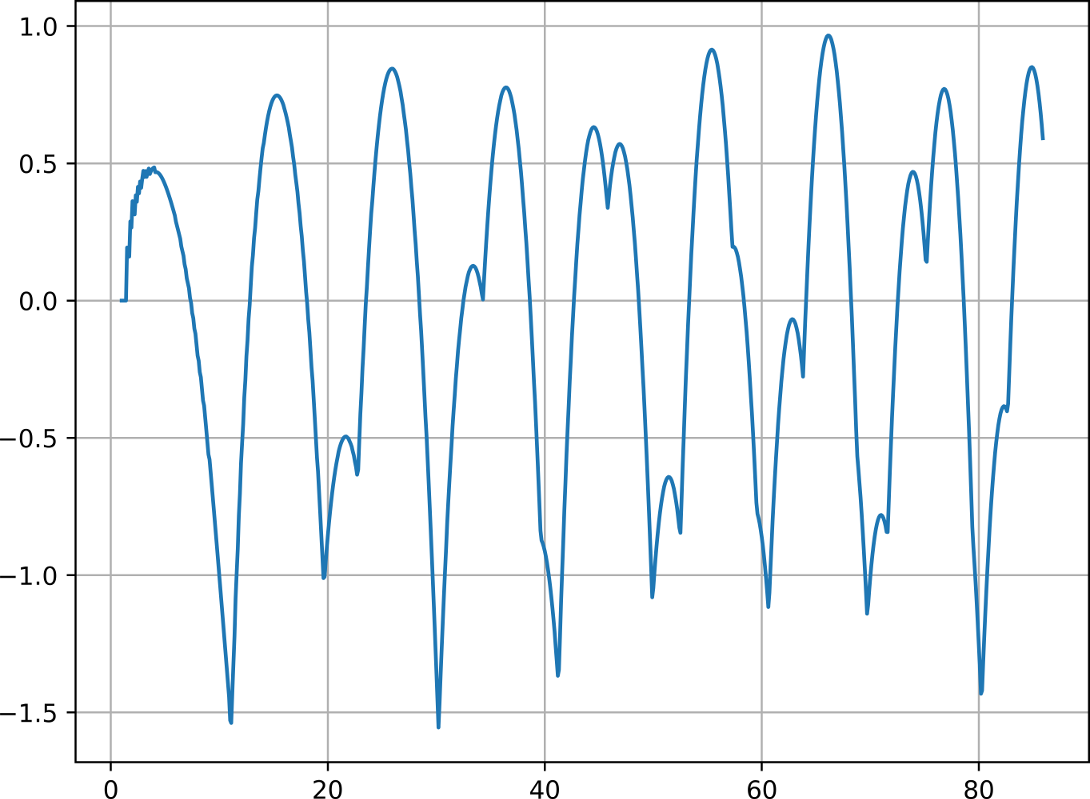}
        \label{fig:Figure2b}
    }
    \caption{Counting function on the tetrahedron}
\end{figure}

\newpage

\section{The Octahedron.}

\begin{figure}[h]
\centering
\includegraphics[width=2in]{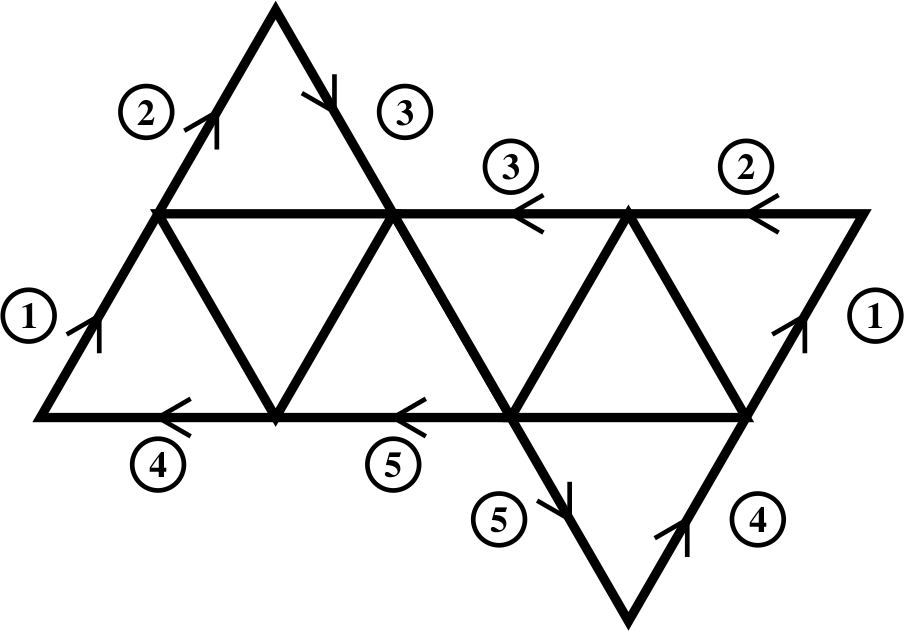}
\caption{Identified edges on the octahedron}
\label{fig:ieoct}
\end{figure}

Figure \ref{fig:ieoct} shows a diagram of the octahedron as a planar region with identified edges.

\begin{figure}[h]
    \centering
    \subfigure[For eigenvalue $\lambda_{57}=16$ (nonsingular)]
    {
        \includegraphics[width=2.2in]{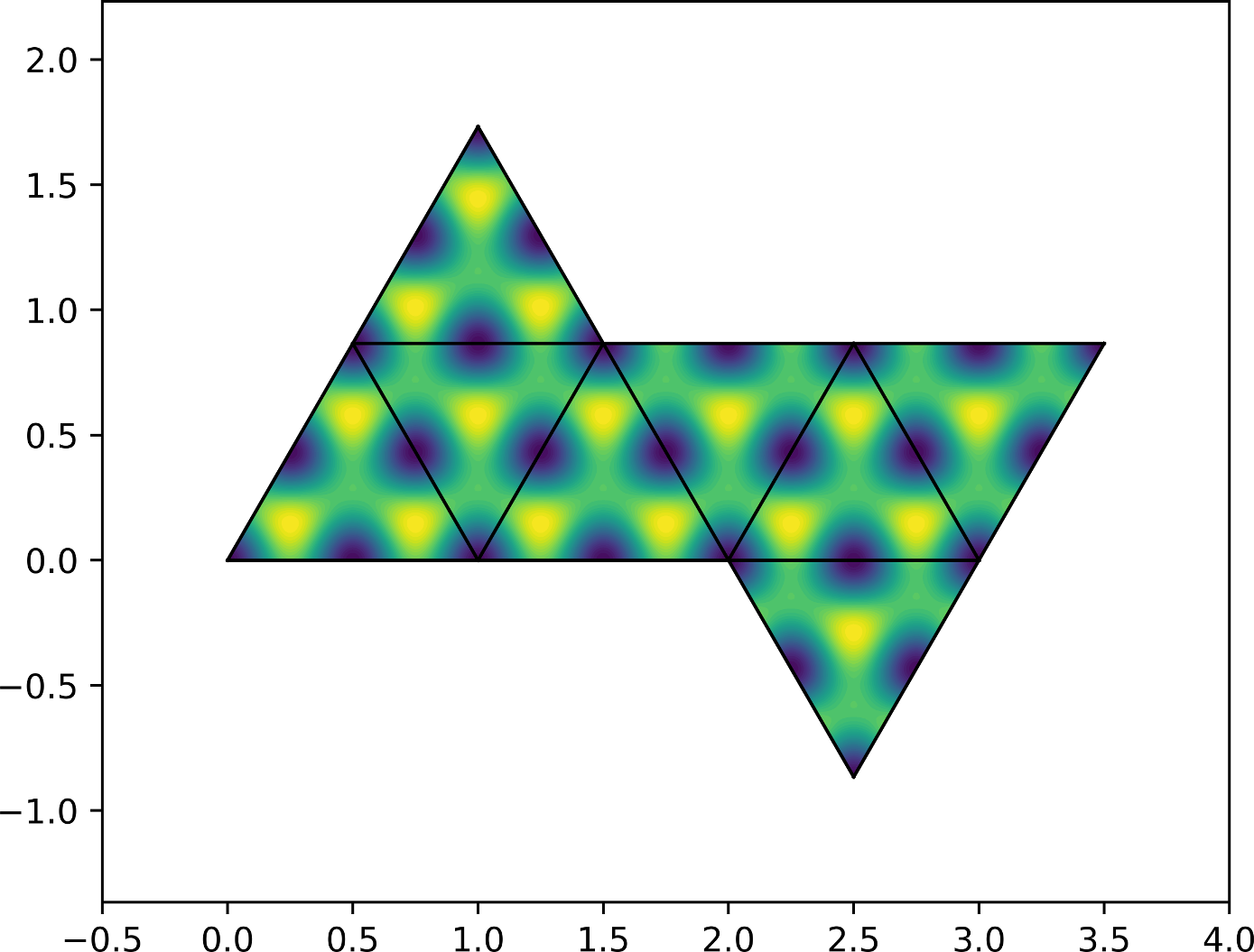}
    }
    \hfill
    \subfigure[For eigenvalue $\lambda_{58}=16.73$ (singular)]
    {
        \includegraphics[width=2.2in]{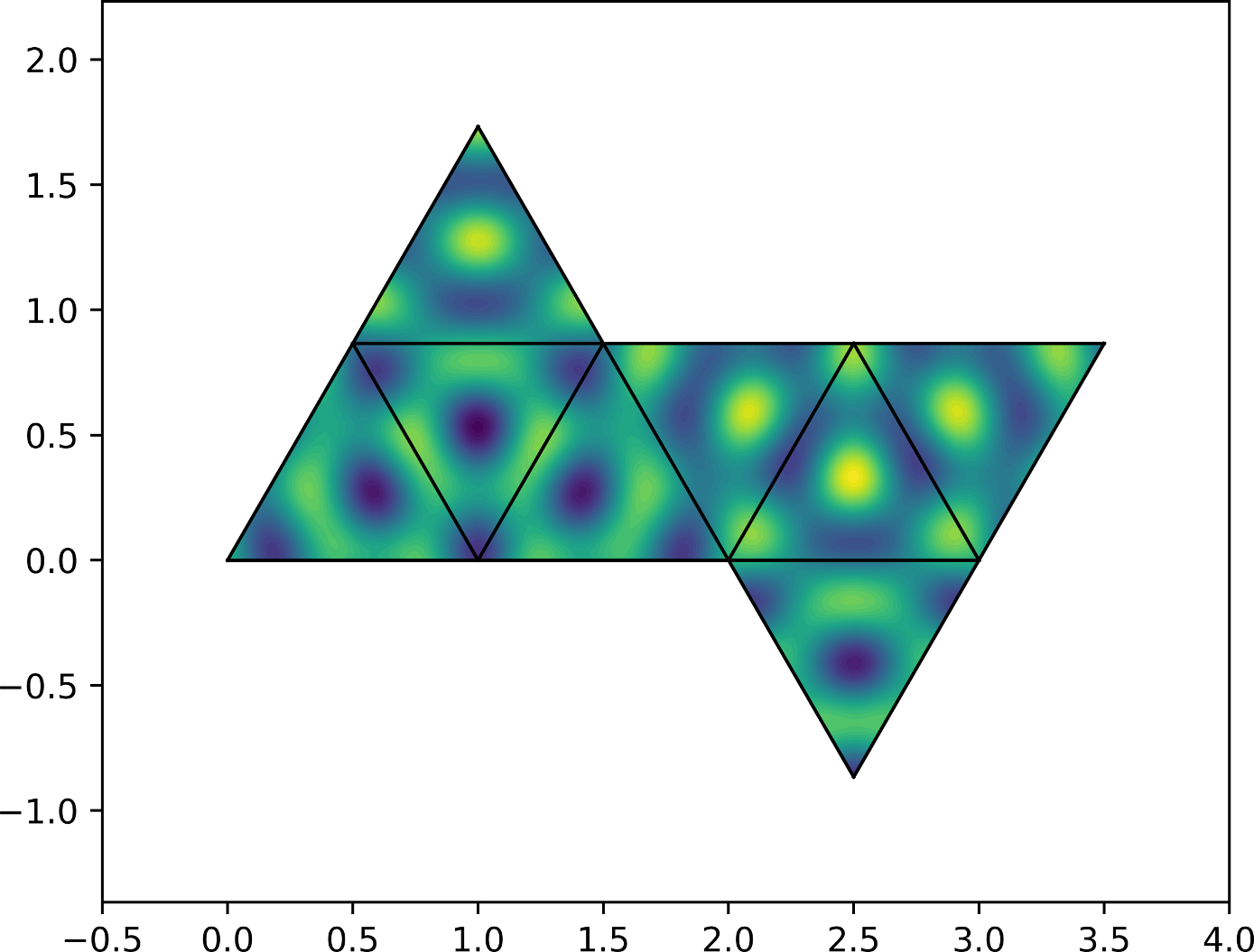}
    }
    \caption{Eigenfunctions on the octahedron}
    \label{fig:pic1oct}
\end{figure}

\begin{figure}[h]
    \centering
    \subfigure[For eigenvalue $\lambda_{57}=16$ (nonsingular)]
    {
        \includegraphics[width=2.2in]{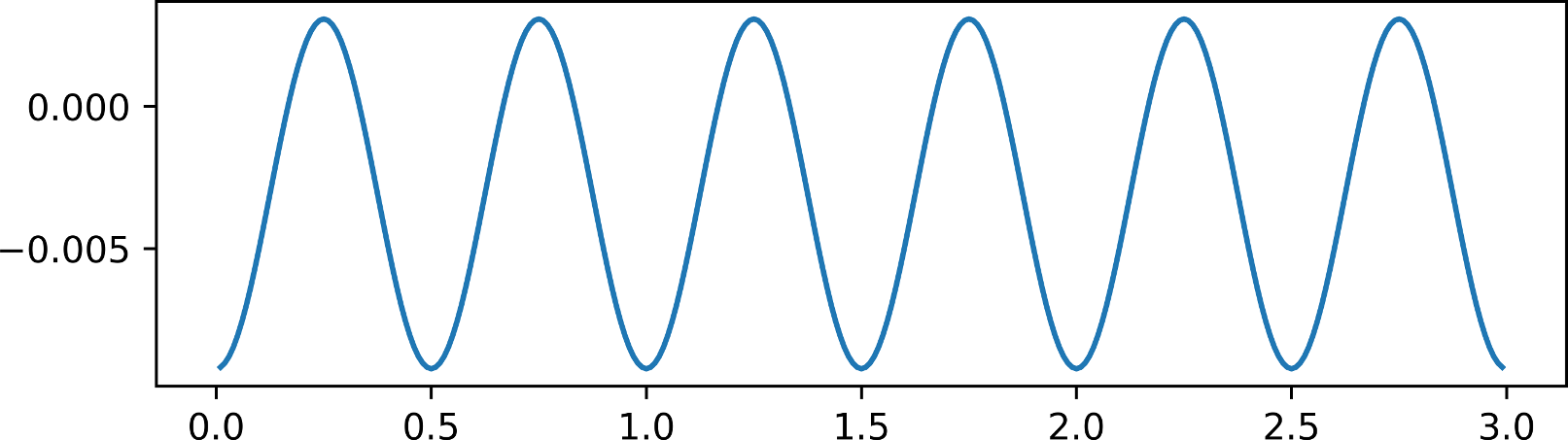}
    }
    \hfill
    \subfigure[For eigenvalue $\lambda_{58}=16.73$ (singular)]
    {
        \includegraphics[width=2.2in]{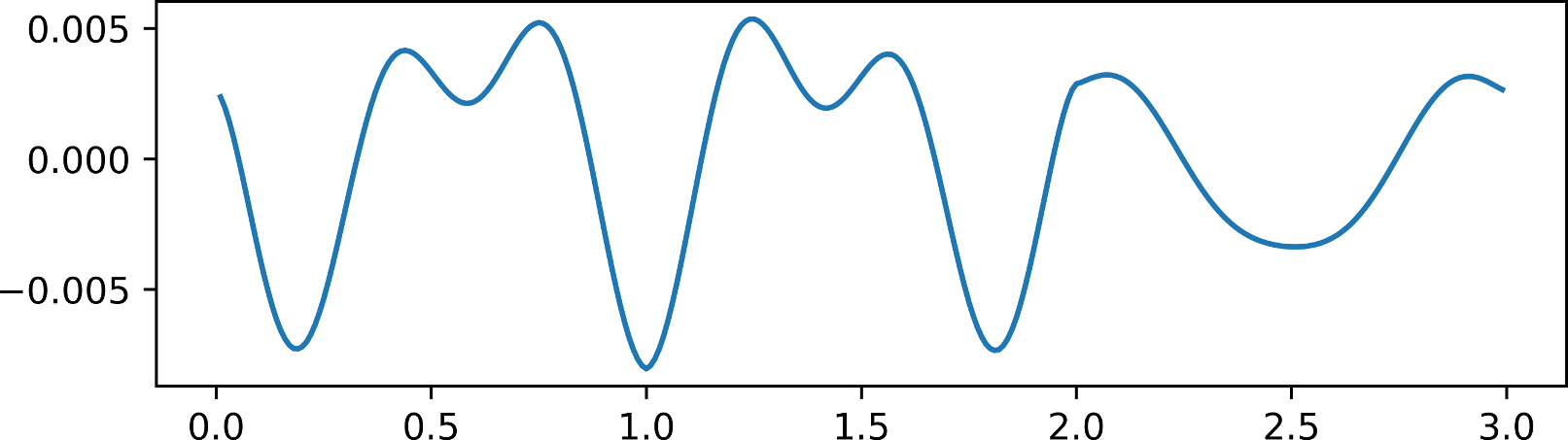}
    }
    \caption{Restriction of eigenfunctions on the octahedron to $y=0$}
    \label{fig:pic2oct}
\end{figure}

We show an example of a nonsingular eigenfunction (Figure \ref{fig:pic1oct}a) and an example of a singular eigenfunction (Figure \ref{fig:pic1oct}b), as well as the respective graph of the restriction to line segments passing through vertices (Figure \ref{fig:pic2oct}).

In Table $2$ we show the beginning of the spectrum of the Laplacian on the surface of the regular octahedron, with the eigenvalues normalized by dividing by $\frac{4}{3}\pi^2$. We extrapolate and will continue using the extrapolated data. After the trivial eigenvalue $0$, all eigenvalues have multiplicities $2$, $3$ or $4$. We notice that all the tetrahedron eigenvalues corresponding to a $1^+$ and $1^-$ eigenfunctions occur -  they are integers and bolded in the table. The multiplicity of the octahedron eigenvalues is always twice the number of the tetrahedron $1^+$ and $1^-$ eigenvalues. We also have eigenvalues $\frac{N}{3}$, where $N$ is an integer eigenvalue - these are italic in the table. The multiplicity of the $\frac{N}{3}$-eigenvalues is the same as the $N$-eigenvalues. We will explain all these empirical observations. We should also mention that the explanations were discovered after we made the observations: experimental mathematics in action!

\vspace{0.15cm}

\begin{table}[h] 
\centering
\makebox[0pt][c]{\parbox{\textwidth}{
\begin{minipage}[b]{0.32\hsize}
\centering
\begin{tabular}{| c | c |}
\hline
\# & Eigenvalue \\ 	
\hline 
  \textbf{1} & \textbf{0} \\
  2 & 0.54376 \\ 
  3 & 0.54376 \\
  4 & 0.54376 \\
  \textit{5} & \textit{1.33342} \\
  \textit{6} & \textit{1.33342} \\
  7 & 1.89224 \\
  8 & 1.89224 \\
  9 & 2.84941 \\
  10 & 2.84941 \\
  11 & 2.84941 \\
  12 & 3.62006 \\
  13 & 3.62006 \\
  14 & 3.62006 \\
  \textbf{15} & \textbf{4.00080} \\
  \textbf{16} & \textbf{4.00080} \\
  \textit{17} & \textit{5.33476} \\
  \textit{18} & \textit{5.33476} \\
  19 & 5.45089 \\
  20 & 5.45089 \\
\hline 
\end{tabular}
\end{minipage}
\hfill
\begin{minipage}[b]{0.32\hsize}
\centering
\begin{tabular}{| c | c |}
\hline
\# & Eigenvalue \\ 	
\hline 
  21 & 5.45089 \\
  22 & 5.45089 \\ 
  23 & 6.37226 \\
  24 & 6.37226 \\
  25 & 6.37226 \\
  26 & 6.84597 \\
  27 & 6.84597 \\
  28 & 6.84597 \\
  29 & 8.38948 \\
  30 & 8.38948 \\
  31 & 8.38948 \\
  32 & 9.18907 \\
  33 & 9.18907 \\
  34 & 9.18907 \\
  \textit{35} & \textit{9.33771} \\
  \textit{36} & \textit{9.33771} \\
  \textit{37} & \textit{9.33771} \\
  \textit{38} & \textit{9.33771} \\
  39 & 10.67867 \\
  40 & 10.67867 \\
\hline 
\end{tabular}
\end{minipage}
\hfill
\begin{minipage}[b]{0.32\hsize}
\centering
\begin{tabular}{| c | c |}
\hline
\# & Eigenvalue \\ 	
\hline 
  41 & 10.67867\\
  \textbf{42} & \textbf{12.00723} \\ 
  \textbf{43} & \textbf{12.00723} \\
  44 & 12.83710 \\
  45 & 12.83710 \\
  46 & 12.83710 \\
  47 & 12.86814 \\
  48 & 12.86814 \\
  49 & 12.86814 \\
  50 & 12.90939 \\
  51 & 12.90939 \\
  52 & 12.90939 \\
  53 & 14.41173 \\
  54 & 14.41173 \\
  55 & 14.41173 \\
  \textbf{56} & \textbf{16.01286} \\
  \textbf{57} & \textbf{16.01286} \\
  58 & 16.72998 \\
  59 & 16.72998 \\
  60 & 16.72998 \\
\hline 
\end{tabular}
\end{minipage}
}}
\vspace{0.5cm}
\caption{Normalized eigenvalues on the octahedron, Res. 128}
\end{table}

\vspace{-0.5cm}

The symmetry group of the octahedron is the direct product of $S_4$ and $\mathbb{Z}_2$. The $S_4$ is the action of all rotations, and they permute pairs of opposite faces. The $\mathbb{Z}_2$ action is just the central reflection $(x_1,x_2,x_3) \to (-x_1,-x_2,-x_3)$ that switches all pairs of opposite vertices. Note that this is a $48$ element group, with twice as many irreducible representations as $S_4$.

\begin{figure}[h]
    \centering
    \subfigure[In-face]
    {
        \includegraphics[width=1.6in]{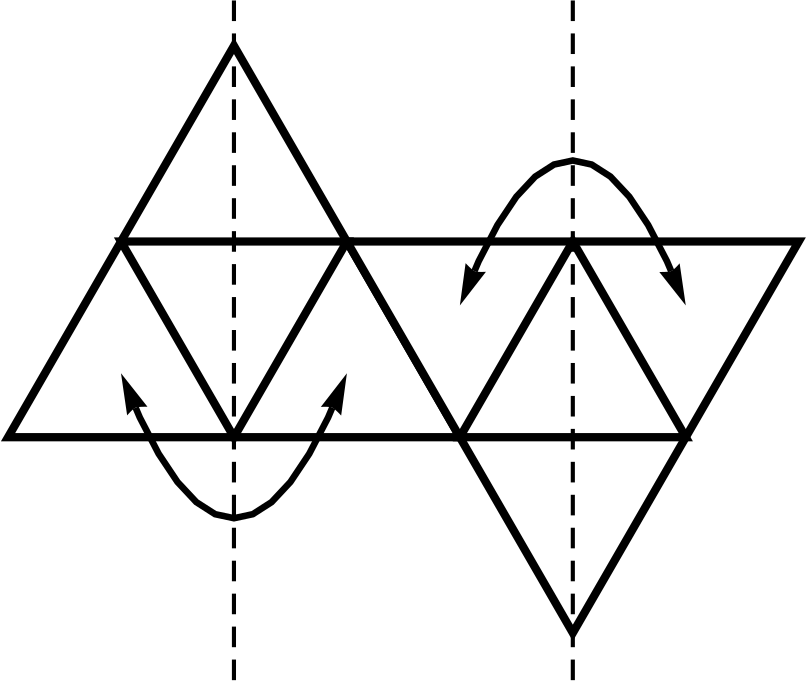}
    }
    \hfill
    \subfigure[Face-to-face]
    {
        \includegraphics[width=1.6in]{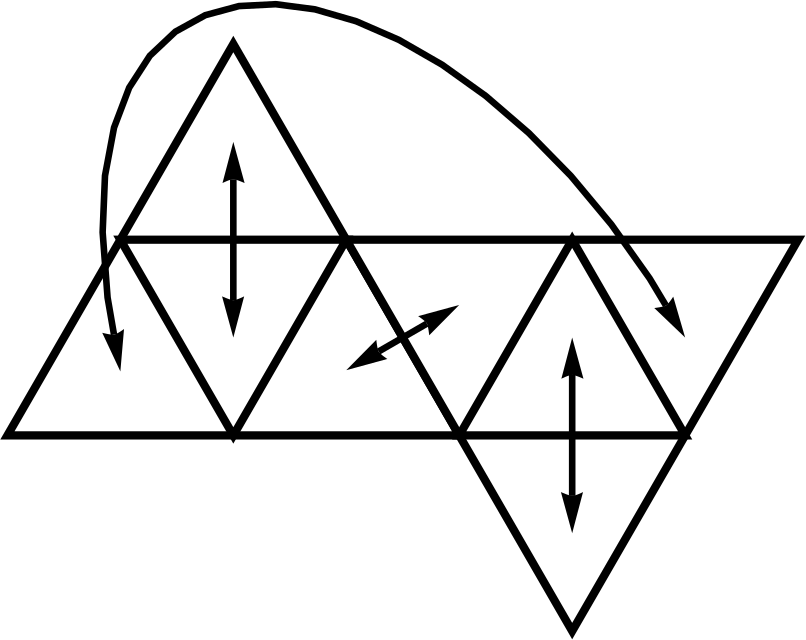}
    }
    \caption{Types of reflections on the octahedron}
    \label{fig:refloct}
\end{figure}

\begin{figure}[h]
    \centering
    \subfigure[$1^{+-}$]
    {
        \includegraphics[width=1.6in]{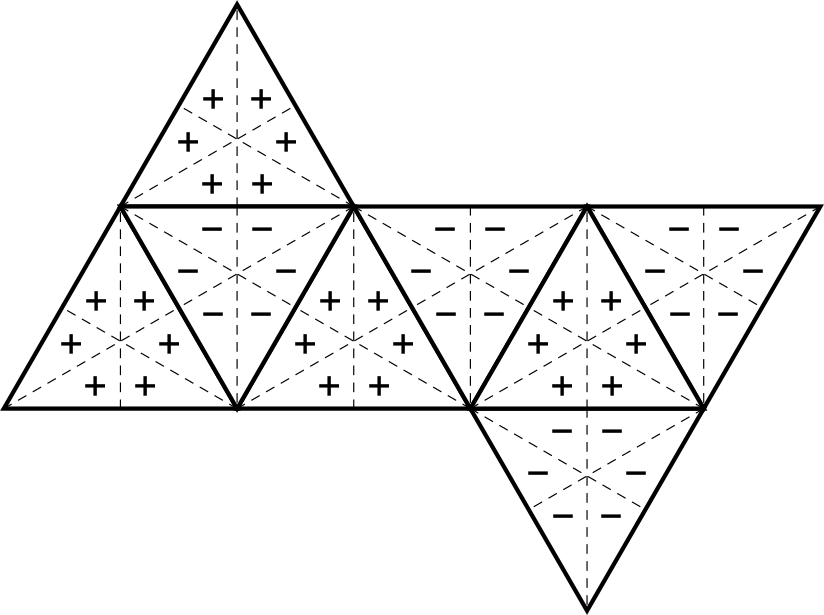}
    }
    \hfill
    \subfigure[$1^{-+}$]
    {
        \includegraphics[width=1.6in]{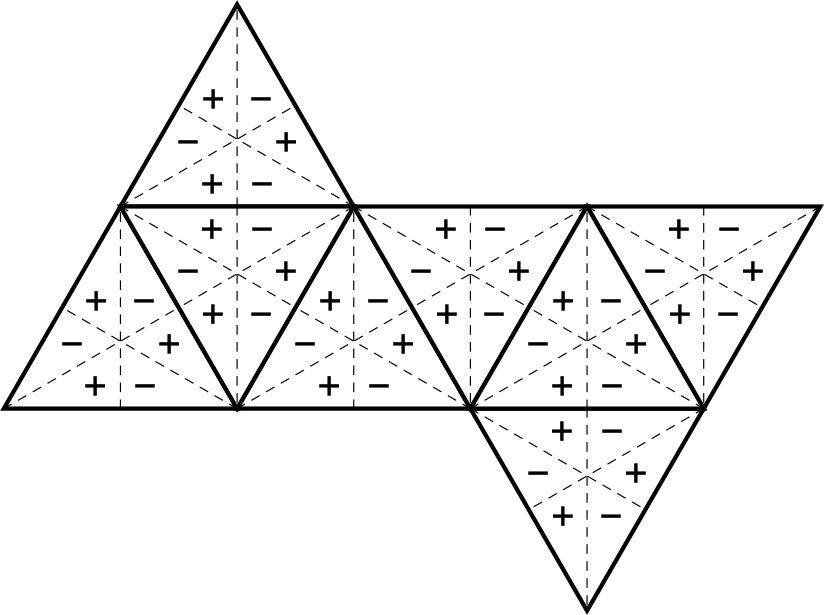}
    }
    \caption{Symmetry properties of eigenfunctions on the octahedron}
    \label{fig:typesoct}
\end{figure}

There are two distinct types of reflections: \textit{in-face} reflections that reflect in the bisectors of $4$ different faces and permute in pairs the remaining faces; and \textit{face-to-face} reflections that reflect four pairs of adjacent faces in their common edges (see Figure \ref{fig:refloct}). In $3$-space the reflection plane passes through two opposite vertices and either slices four faces in half or passes through four edges. So there are four $1$-dimensional representations that we denote by $1^{\pm \pm}$, the first $\pm$ choice indicating symmetry or skew-symmetry with respect to in-face reflections, and the second $\pm$ choice indicating the same with respect to face-to-face reflections. The $1^{++}$ eigenfunctions have symmetry with respect to all reflections, while the $1^{--}$ eigenfunctions have skew-symmetry. The $1^{+-}$ and $1^{-+}$ symmetry properties are illustrated in Figure \ref{fig:typesoct}. As in the case of the tetrahedron, a fundamental domain may be taken to be a one-sixth face.

\vspace{-0.35cm}

\begin{figure}[h]
    \centering
    \subfigure[Eigenfunction for eigenvalue $\lambda_{85}=48$ on the tetrahedron]
    {
        \includegraphics[width=2.16in]{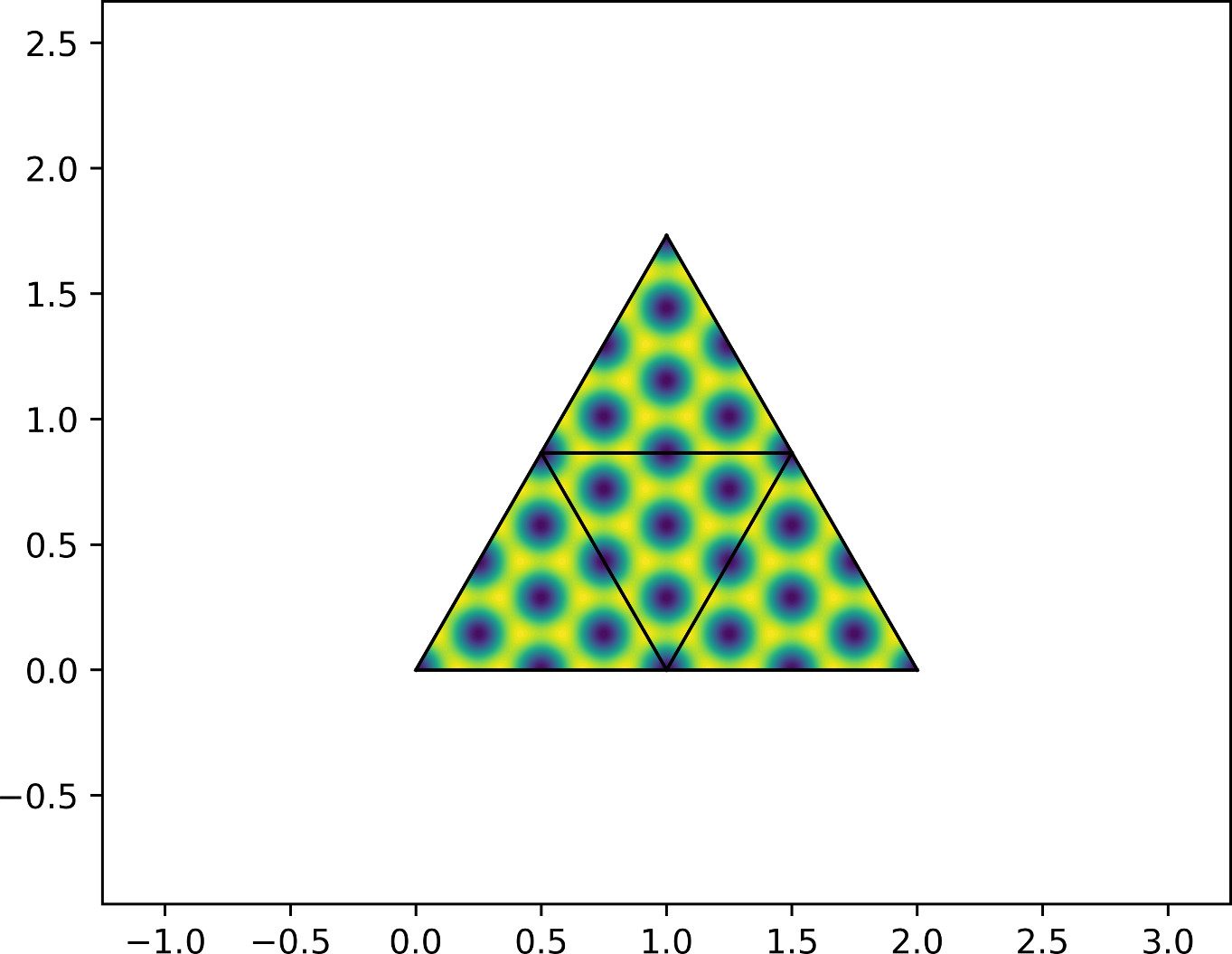}
    }
    \hfill
    \subfigure[Eigenfunction for eigenvalue $\lambda_{171}=48$ on the octahedron]
    {
        \includegraphics[width=2.2in]{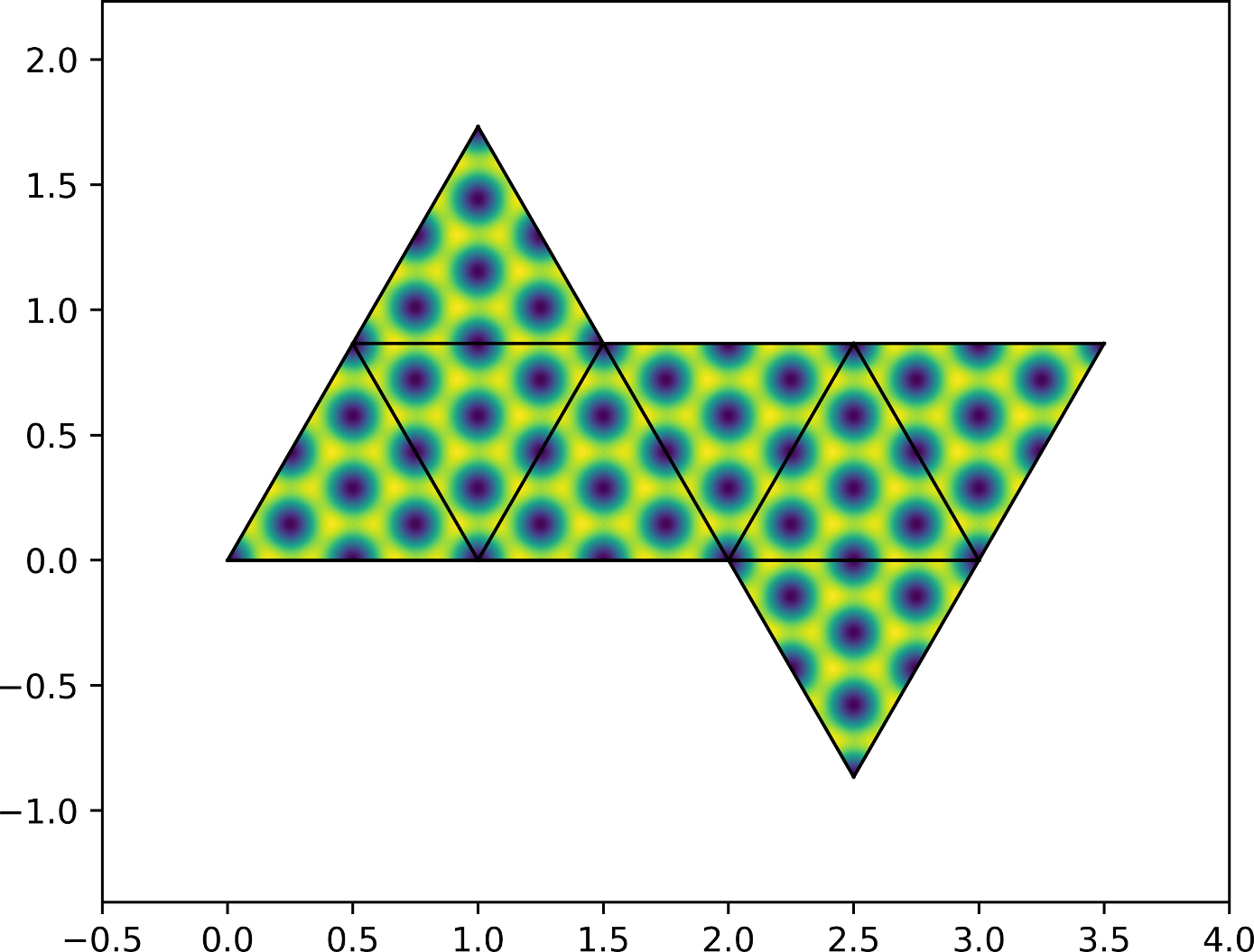}
    }
    \caption{Tetrahedron-type eigenfunctions}
    \label{fig:tt}
\end{figure}

It is clear from inspection that all four types are well-defined on the small two-faced torus. This may also be observed as a consequence of the fact that the torus translations are products of two face-to-face reflections. So for the same reason as for the tetrahedron, all four eigenfunction types must be linear combinations of (\ref{eq:torfuncs}) with both $k$ and $j$ even. In particular, $u_{\pm}$ given by (\ref{eq:uplusmin}), (\ref{eq:uplus}), and (\ref{eq:uplus2}) can also be regarded as $1^{++}$ and $1^{--}$ eigenfunctions on the octahedron. We will refer to these as \textit{tetrahedron-type} eigenfunctions, as they are identical on faces of both polyhedra (Figure \ref{fig:tt}).

\vspace{-0.1cm}

\newpage

\begin{figure}[h]
    \centering
    \subfigure[$1^{+-}$ \hspace{.1cm} Face-to-face reflection axes]
    {
        \includegraphics[width=1.8in]{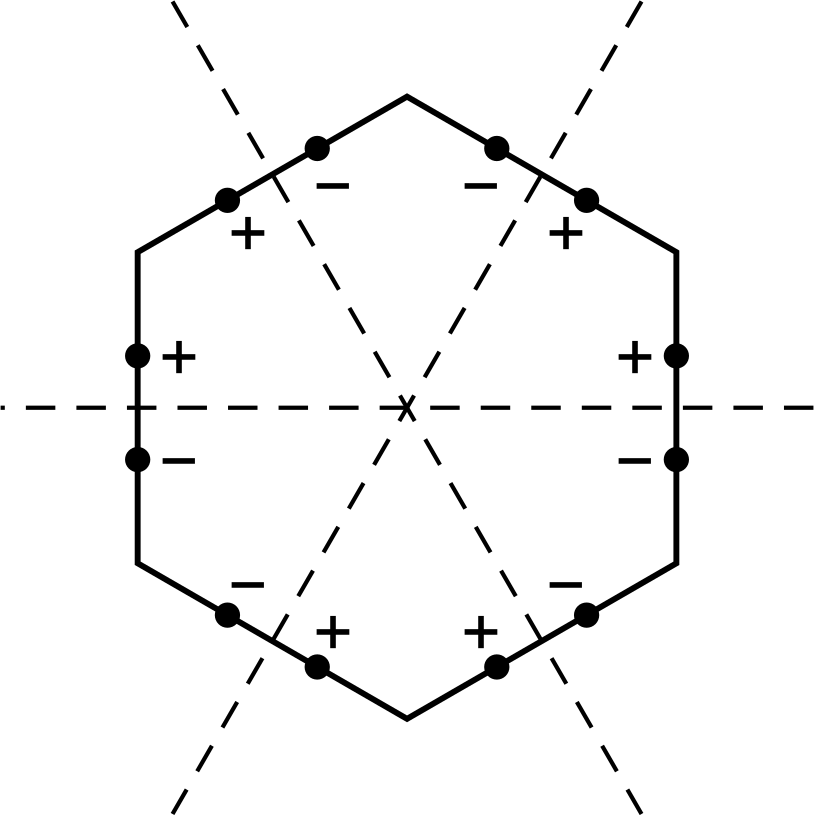}
    }
    \subfigure[$1^{-+}$ \hspace{.1cm}  In-face reflection axes]
    {
        \includegraphics[width=1.8in]{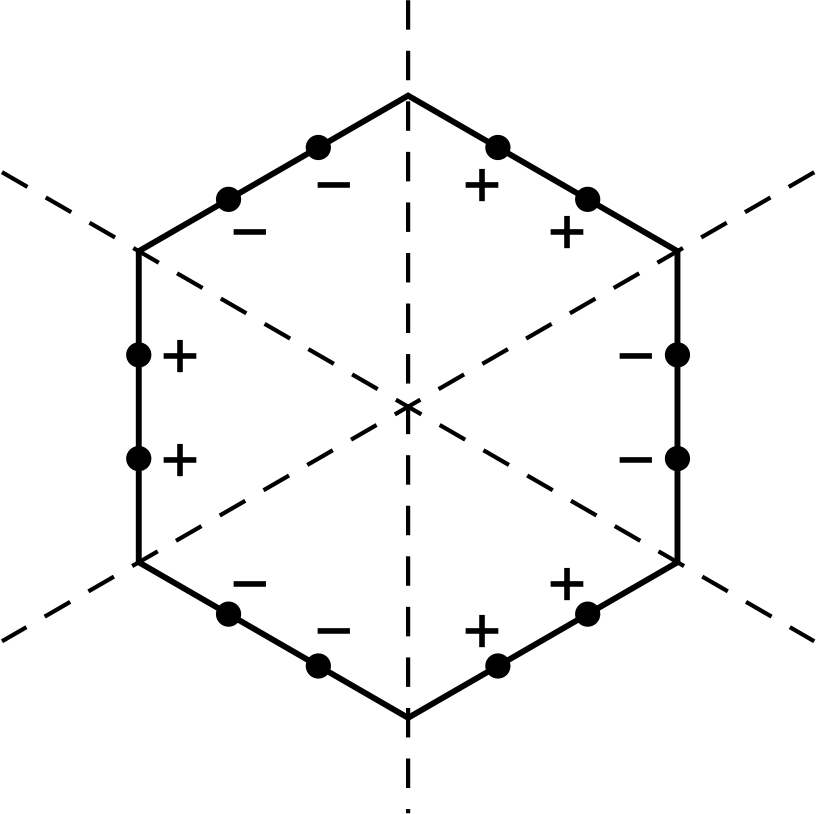}
    }
    \caption{Distribution of signs}
    \label{fig:distoct}
\end{figure}

To understand the $1^{+-}$ and $1^{-+}$ eigenfunctions we need to return to the diagram of a generic orbit in Figure $12$. The $+$ and $-$ labels there correspond to $1^{--}$ eigenfunctions. Note that from Figure \ref{fig:typesoct} we see that, under the half turn symmetry of the torus, the $1^{+-}$ and $1^{-+}$ eigenfunctions are skew-symmetric. That means that the cosines in (\ref{eq:torfuncs}) must be replaced by sines. This is also clear from the fact that opposite lattice points have opposite signs in Figure \ref{fig:distoct}, which shows the distribution of signs for the two types and also shows the reflection axes corresponding to the two types of reflections.

\noindent From this we may read off the explicit formulas (again $k>j>0$, both even):
\begin{equation*}
\begin{aligned}
u_{+-} &= \sin(2\pi x \cdot (k\overset{\rightarrow}{u} + j\overset{\rightarrow}{v}))-\sin(2\pi x \cdot (j\overset{\rightarrow}{u} + k \overset{\rightarrow}{v})) \\
       &+ \sin(2\pi x \cdot ((k+j)\overset{\rightarrow}{u} - j \overset{\rightarrow}{v})) - \sin(2\pi x \cdot ((k+j)\overset{\rightarrow}{u} - k\overset{\rightarrow}{v})) \\
       &+ \sin(2\pi x \cdot (j\overset{\rightarrow}{u} - (j+k) \overset{\rightarrow}{v})) - \sin(2\pi x \cdot (k\overset{\rightarrow}{u} - (k+j)\overset{\rightarrow}{v}))\\
\end{aligned}    
\end{equation*}    
       
\noindent and

\begin{equation*}
\begin{aligned}
u_{-+} &= \sin(2\pi x \cdot (k\overset{\rightarrow}{u} + j\overset{\rightarrow}{v}))-\sin(2\pi x \cdot (j\overset{\rightarrow}{u} + k \overset{\rightarrow}{v})) \\
       &- \sin(2\pi x \cdot ((k+j)\overset{\rightarrow}{u} - j \overset{\rightarrow}{v})) - \sin(2\pi x \cdot ((k+j)\overset{\rightarrow}{u} - k\overset{\rightarrow}{v})) \\
       &+ \sin(2\pi x \cdot (j\overset{\rightarrow}{u} - (j+k) \overset{\rightarrow}{v})) + \sin(2\pi x \cdot (k\overset{\rightarrow}{u} - (k+j)\overset{\rightarrow}{v}))
\end{aligned}
\end{equation*}

\vspace{0.3cm}

\noindent For nongeneric orbits only one of these survives and collapses to three terms:
\begin{equation}
\label{eq:three3}
u_{+-} = \sin(2\pi x \cdot (k\overset{\rightarrow}{u}))- \sin(2\pi x \cdot (k\overset{\rightarrow}{v})) - \sin(2\pi x \cdot (k\overset{\rightarrow}{u}-k\overset{\rightarrow}{v}))
\end{equation}

\noindent This is a type $1$ nongeneric orbit with $j=0$.
\begin{equation}
\label{eq:three4}
u_{-+} = \sin(2\pi x \cdot (k\overset{\rightarrow}{u} + k \overset{\rightarrow}{v})) - \sin(2\pi x \cdot (2k\overset{\rightarrow}{u} - k\overset{\rightarrow}{v})) + \sin(2\pi x \cdot (k\overset{\rightarrow}{u} - 2k\overset{\rightarrow}{v}))
\end{equation}

\noindent This is a type $2$ nongeneric orbit with $j=k$.

Note that these eigenfunctions restricted to a face are distinct from all the tetrahedron eigenfunctions on a face because they have the skew-symmetry with respect to the half-turn. Nevertheless they share the same eigenvalue $\frac{4\pi^2}{3}(j^2 + k^2 +jk)$. This explains the doubling of the mutliplicities: eigenspaces corresponding to generic orbits split into four distinct symmetry types, while those corresponding to nongeneric orbits split into two distinct symmetry types. And of course all these eigenfunctions are nonsingular.

But these are not the only nonsingular eigenfunctions. There is another family that is derived from these by a process of so called \emph{enlargement}. Note that the fundamental domain is a $30^\circ-60^\circ-90^\circ$ triangle that is similar to a half-face, and the contraction ratio is $\sqrt{3}$. If we take $u$ to be any of $u_+,u_-,u_{+-},u_{-+}$ and compose with this similarity $S:$ half-face $\to$ fundamental domain, then $u \circ S$ is an eigenfunction of the Laplacian on the half-face with eigenvalue multiplied by $\frac{1}{3}$, and we may extend it to the octahedron by following the same reflection rules. Figure \ref{fig:typesoct} shows the $--$, $-+,$ and $+-$ cases. Note that there are three choices in how to bisect a face into two half-faces, and we make these choices above to be consistent. When we make the identifications as shown in Figure $18$, the reflection rules are preserved. In Table $2$ we also identified the eigenvalues of the form $\frac{N}{3}$ and saw that the multiplicity agrees with the multiplicity of $N$. When the integer $N$ is divisible by $3$ (first example $N=12$), the eigenspace created by enlargement coincides with the one previously constructed. Otherwise, the $\frac{N}{3}$ eigenspace does not split into $1\pm \pm$ type eigenfunctions.

We again show the counting function $N(t)$, the remainder $D(t)$, and the averages $A(t)$ and $g(t)$ in Figure 24. The constant $\frac{\sqrt{3}}{2\pi}$ comes from the Weyl term and the constant $\frac{5}{12}=6\cdot\frac{5}{72}$ is the value conjectured in [5].

\begin{figure}[h]
    \centering
    \subfigure[The counting function $N(t)$ (blue) and $\left(\frac{\sqrt{3}}{2\pi}t + \frac{5}{12}\right)$ (orange)]
    {
        \includegraphics[width=2.1in]{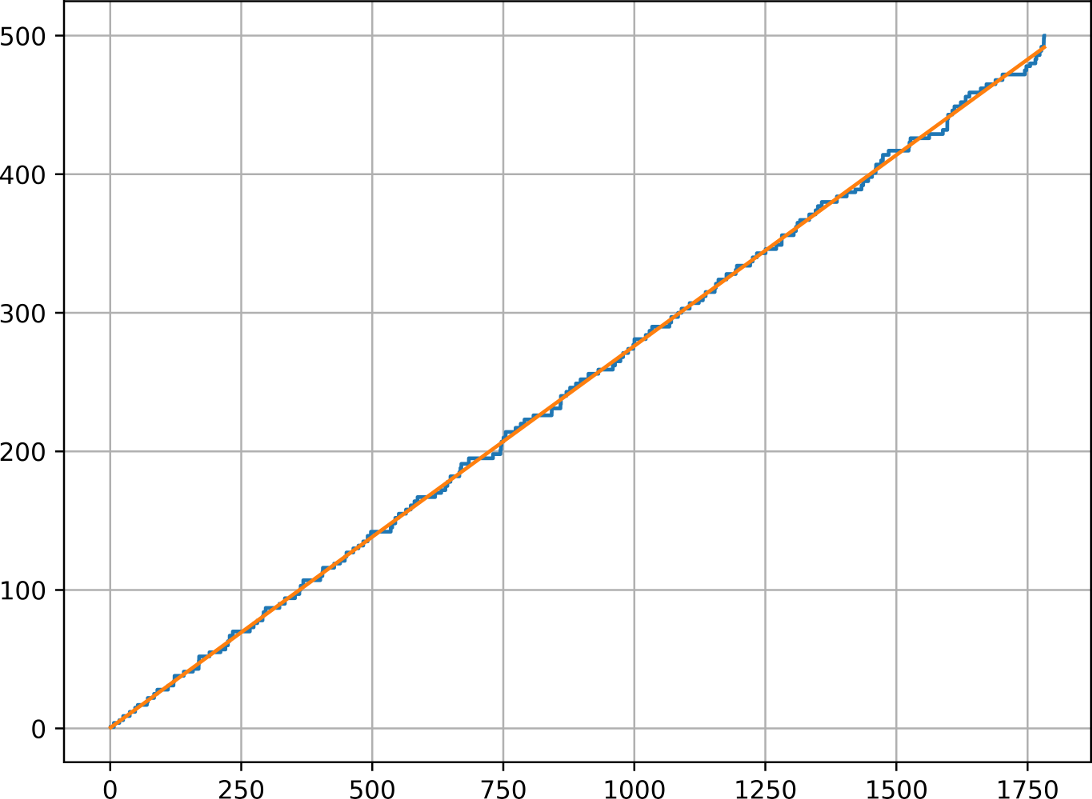}
        \label{fig:Figure1a}
    }
    \hfill
    \subfigure[The difference $D(t)$]
    {
        \includegraphics[width=2.1in]{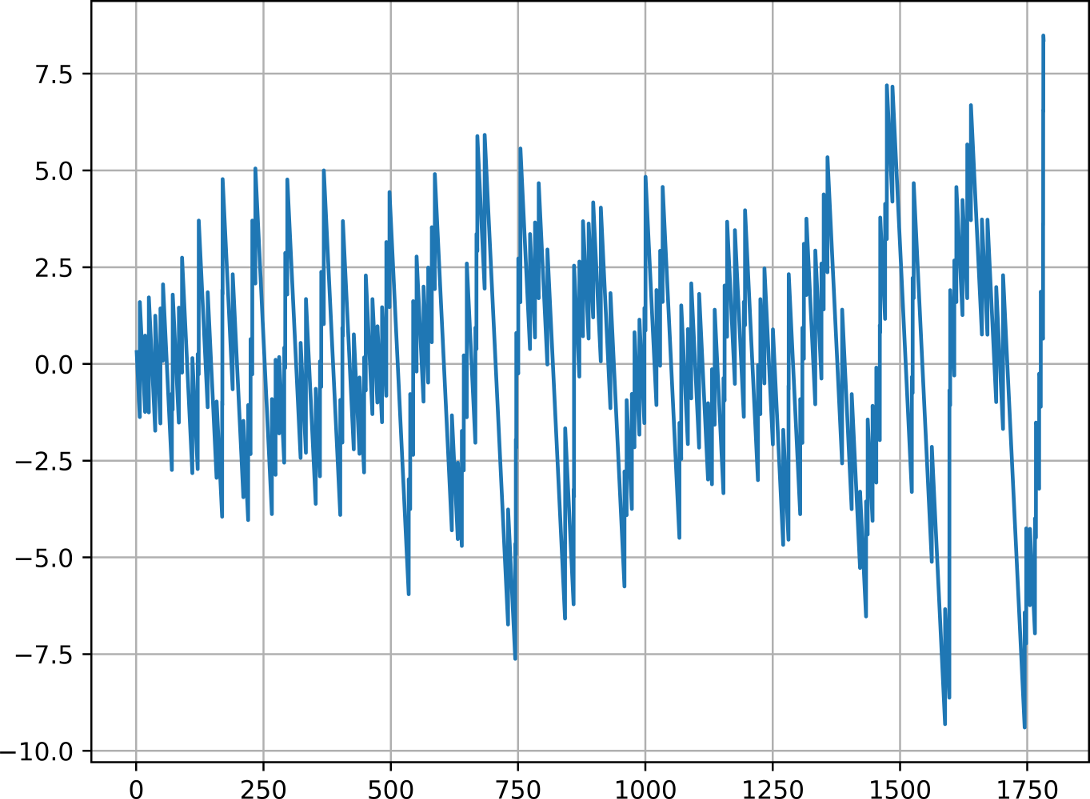}
        \label{fig:Figure2b}
    }
    \subfigure[The average of the difference $A(t)$]
    {
        \includegraphics[width=2.1in]{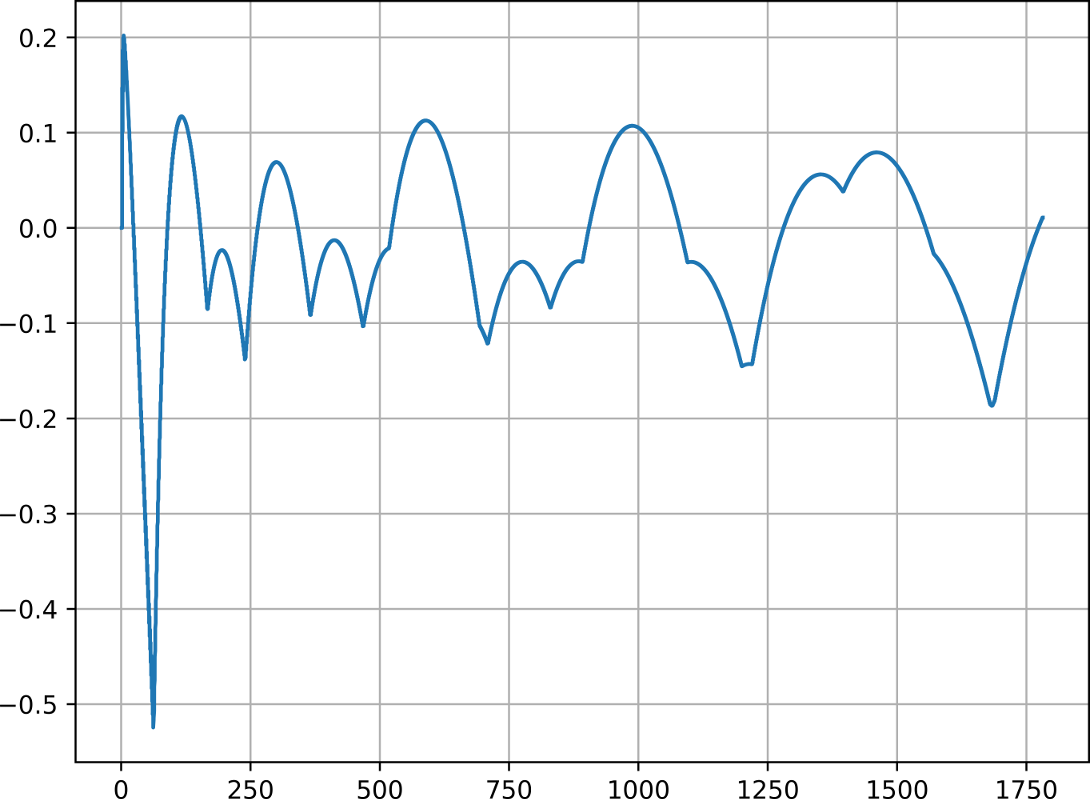}
        \label{fig:Figure1a}
    }
    \hfill
    \subfigure[The average after rescaling $g(t)$]
    {
        \includegraphics[width=2.1in]{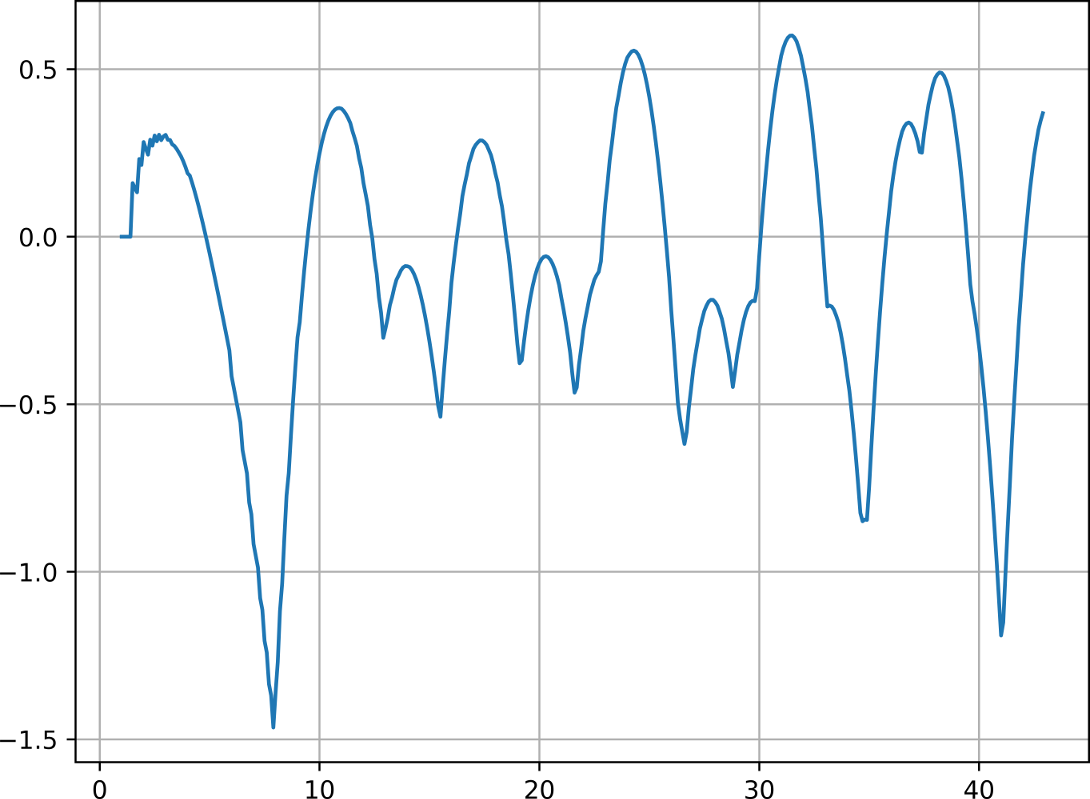}
        \label{fig:Figure2b}
    }
    \caption{Counting function on the octahedron}
\end{figure}

\newpage

\section{The Icosahedron.}

\begin{figure}[h]
    \centering
    \subfigure[For eigenvalue $\lambda_{110}=12$ (nonsingular)]
    {
        \includegraphics[width=2.2in]{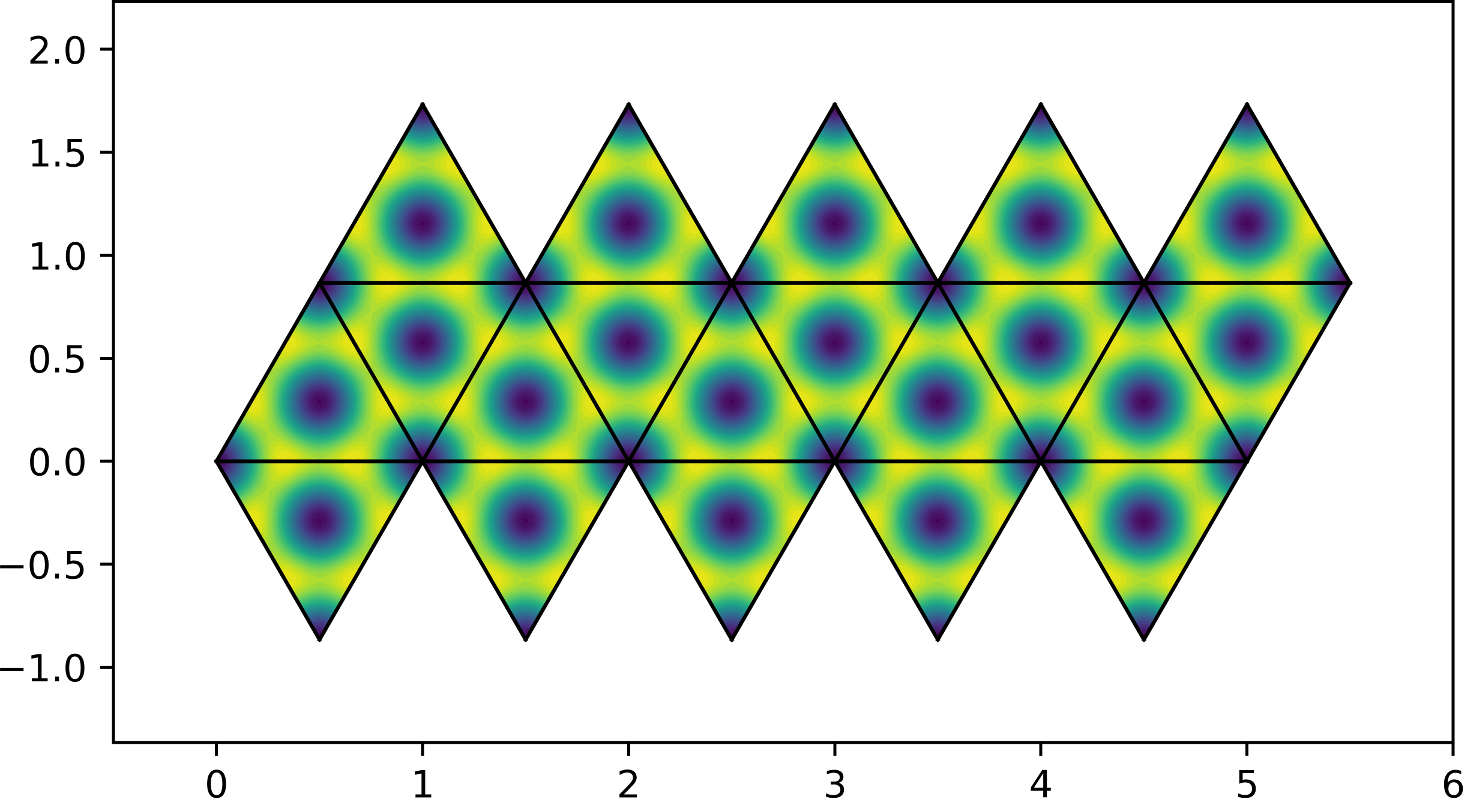}
    }
    \subfigure[For eigenvalue $\lambda_{111}=12.48$ (singular)]
    {
        \includegraphics[width=2.2in]{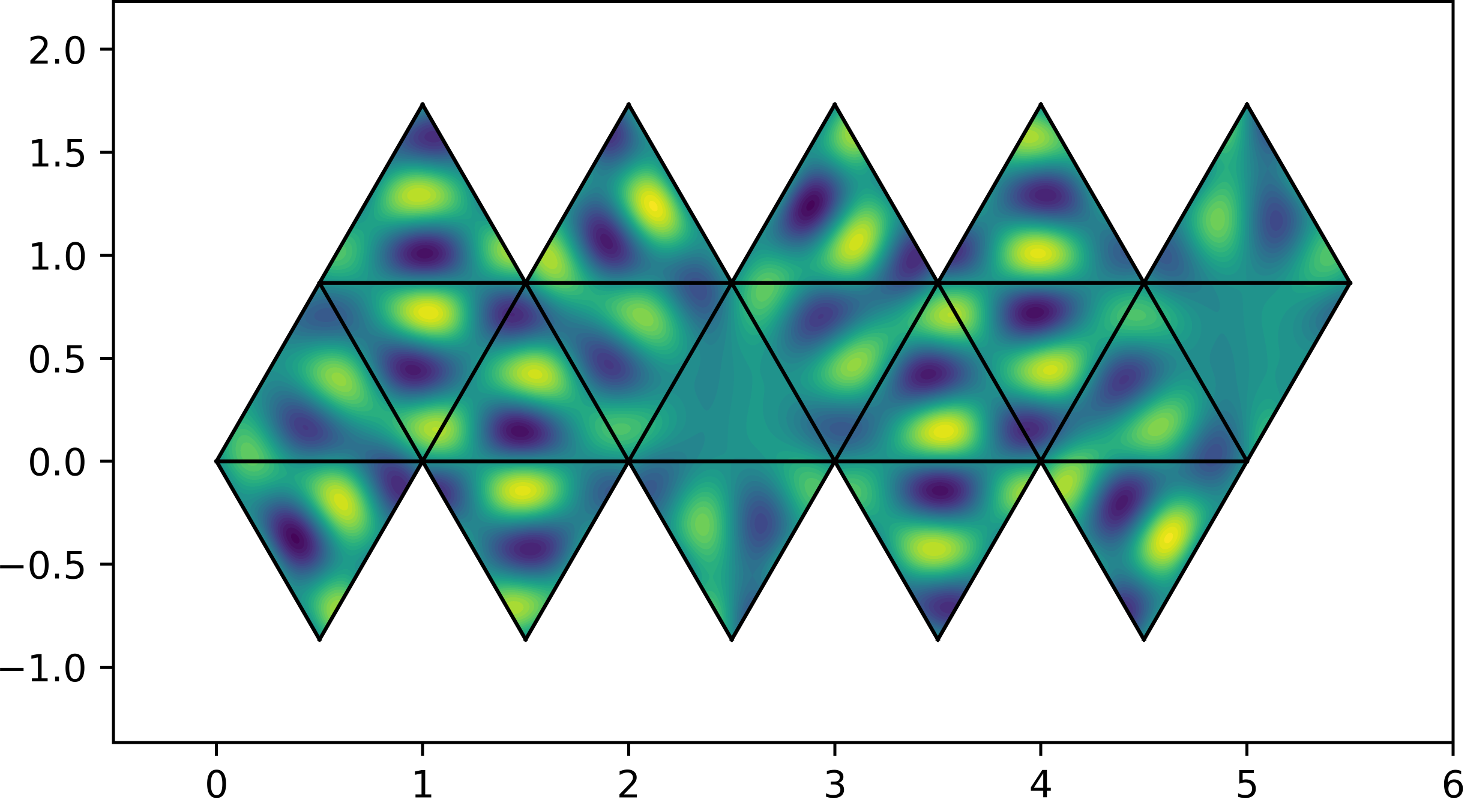}
    }
    \caption{Eigenfunctions on the icosahedron}
    \label{fig:icos1}
\end{figure}

\begin{figure}[h]
    \centering
    \subfigure[For eigenvalue $\lambda_{110}=12$ (nonsingular)]
    {
        \includegraphics[width=2.2in]{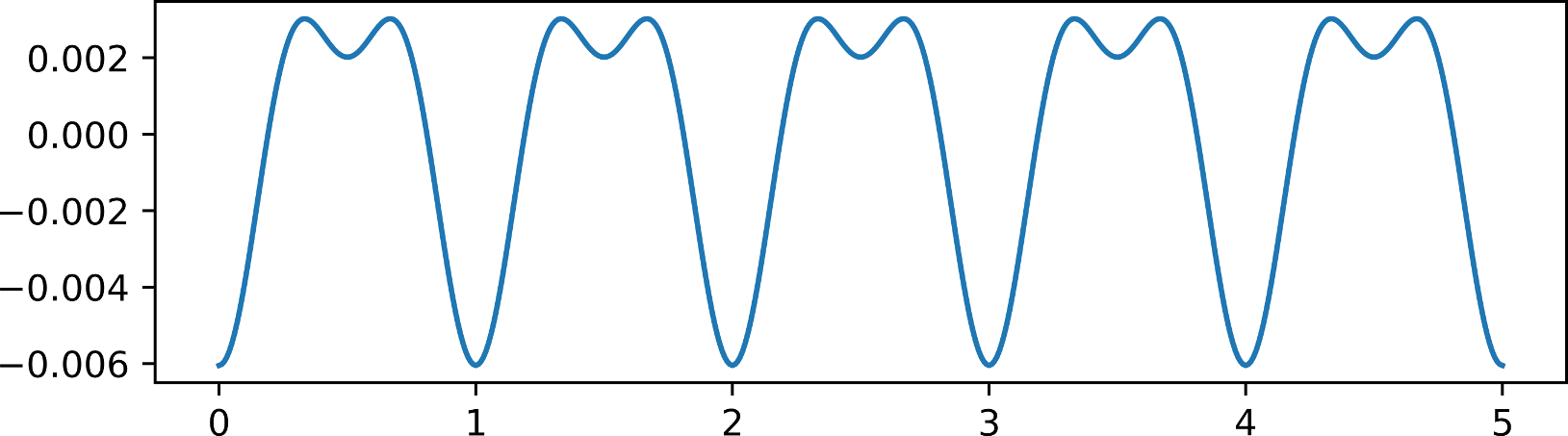}
    }
    \subfigure[For eigenvalue $\lambda_{111}=12.48$ (singular)]
    {
        \includegraphics[width=2.2in]{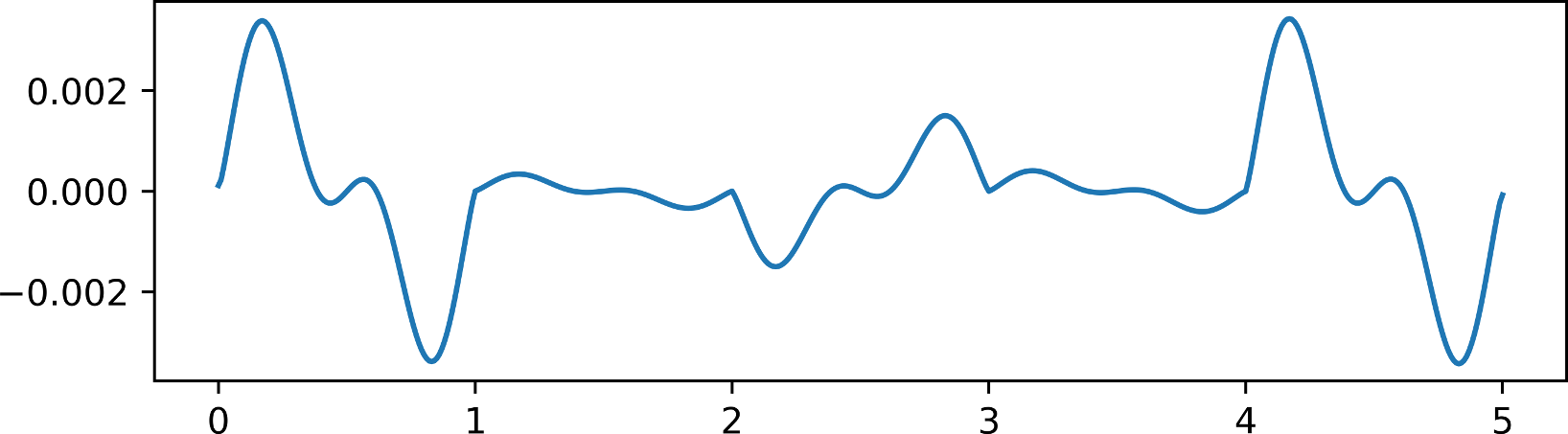}
    }
    \caption{Restriction of eigenfunctions on the icosahedron to $y=0$}
    \label{fig:icos2}
\end{figure}

We start by showing an example of a nonsingular eigenfunction (Figure \ref{fig:icos1}a) and an example of a singular eigenfunction (Figure \ref{fig:icos1}b), as well as the respective graph of the restriction to line segments passing through vertices (Figure \ref{fig:icos2}).

\begin{figure}[h]
\centering
\includegraphics[width=2.5in]{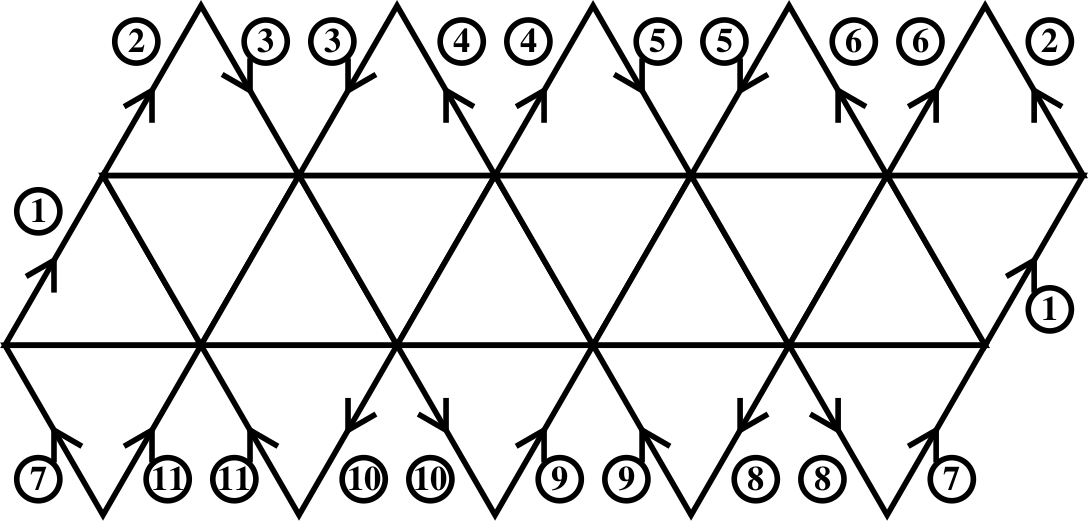}
\caption{Identified edges on the icosahedron}
\label{fig:icosie}
\end{figure}

The surface of the regular icosahedron may be represented as a planar region with identified edges as shown in Figure \ref{fig:icosie}. Each vertex is incident to five faces. Since five is odd, the icosahedron is more closely related to the tetrahedron than the octahedron. The symmetry group is $A_5 \times \mathbb{Z}_2$, where $A_5$ is the alternating subgroup of $S_5$, which acts via rotations. The $\mathbb{Z}_2$ component is generated by $-I$, as in the case of the octahedron. However, since $A_5$ only has the trivial one-dimensional representation, there are only two irreducible one-dimensional representations, yielding the $1^+$ symmetry type of functions symmetric with respect to all reflections, and the $1^-$ type skew-symmetric with respect to all reflections.

As in the case of the tetrahedron, reflections that reflect in-face for some faces will swap other faces by reflecting in their common edge.

The icosahedron will have $1^+$ and $1^-$ eigenfunctions of tetrahedral type. On each face they will be given by (\ref{eq:uplusmin}), (\ref{eq:uplus}) and (\ref{eq:uplus2}), and they will extend to the other faces by even ($+$) or odd ($-$) reflections. These are nonsingular, and they appear to be the only nonsingular eigenfunctions. The enlgargement idea cannot be extended consistently across the entire icosahedron. Table $3$ lists the beginning portion of the spectrum, normalized as before with integer eigenvalues bolded. We extrapolate and will continue using the extrapolated data. The group $A_5$ has five distinct irreducible representations, of dimensions $1$, $3$, $3$, $4$, $5$. In the Table we see multiplicities $3$, $4$ and $5$ for the noninteger values. In this case there appear to be no coincidences (distinct representation types sharing a common eigenvalue).

\vspace{0.1cm}

\begin{table}[ht] 
\centering
\makebox[0pt][c]{\parbox{\textwidth}{
\begin{minipage}[b]{0.32\hsize}
\centering
\begin{tabular}{| c | c |}
\hline
\# & Eigenvalue \\ 	
\hline 
  \textbf{1} & \textbf{0} \\
  2 & 0.22032 \\ 
  3 & 0.22032 \\
  4 & 0.22032 \\
  5 & 0.65895 \\
  6 & 0.65895 \\
  7 & 0.65896 \\
  8 & 0.65896 \\
  9 & 0.65895 \\
  10 & 1.22415 \\
  11 & 1.22415 \\
  12 & 1.22415 \\
  13 & 1.39760 \\
  14 & 1.39761 \\
  15 & 1.39761 \\
  16 & 1.39762 \\
  17 & 2.11275 \\
  18 & 2.11276 \\
  19 & 2.11277 \\
  20 & 2.11277 \\
\hline 
\end{tabular}
\end{minipage}
\hfill
\begin{minipage}[b]{0.32\hsize}
\centering
\begin{tabular}{| c | c |}
\hline
\# & Eigenvalue \\ 	
\hline 
  21 & 2.11277 \\
  22 & 2.32749 \\ 
  23 & 2.32749 \\
  24 & 2.32750 \\
  25 & 2.32750 \\
  26 & 3.05440 \\
  27 & 3.05441 \\
  28 & 3.05442 \\
  29 & 3.40530 \\
  30 & 3.40530 \\
  31 & 3.40530 \\
  32 & 3.40727 \\
  33 & 3.40728 \\
  34 & 3.40730 \\
  35 & 3.40732 \\
  36 & 3.40732 \\
  \textbf{37} & \textbf{4.00080} \\
  38 & 4.56435 \\
  39 & 4.56436 \\
  40 & 4.56439 \\ 
\hline 
\end{tabular}
\end{minipage}
\hfill
\begin{minipage}[b]{0.32\hsize}
\centering
\begin{tabular}{| c | c |}
\hline
\# & Eigenvalue \\ 	
\hline 
  41 & 4.64893 \\
  42 & 4.64893 \\ 
  43 & 4.64893 \\
  44 & 4.64894 \\
  45 & 4.64894 \\
  46 & 4.83217 \\
  47 & 4.83219 \\
  48 & 4.83221 \\
  49 & 4.83222 \\
  50 & 5.69309 \\
  51 & 5.69313 \\
  52 & 5.69313 \\
  53 & 6.19595 \\
  54 & 6.19595 \\
  55 & 6.19595 \\
  56 & 6.27054 \\
  57 & 6.27057 \\
  58 & 6.27060 \\
  59 & 6.27062 \\
  60 & 6.27062 \\
\hline 
\end{tabular}
\end{minipage}
}}
\vspace{0.5cm}
\caption{Normalized eigenvalues on the icosahedron, Res. 128}
\end{table}

\vspace{-0.4cm}

\newpage

\begin{figure}[h]
    \centering
    \subfigure[The counting function $N(t)$ (blue) and $\left(\frac{5\sqrt{3}}{4\pi}t + \frac{11}{30}\right)$ (orange)]
    {
        \includegraphics[width=2.1in]{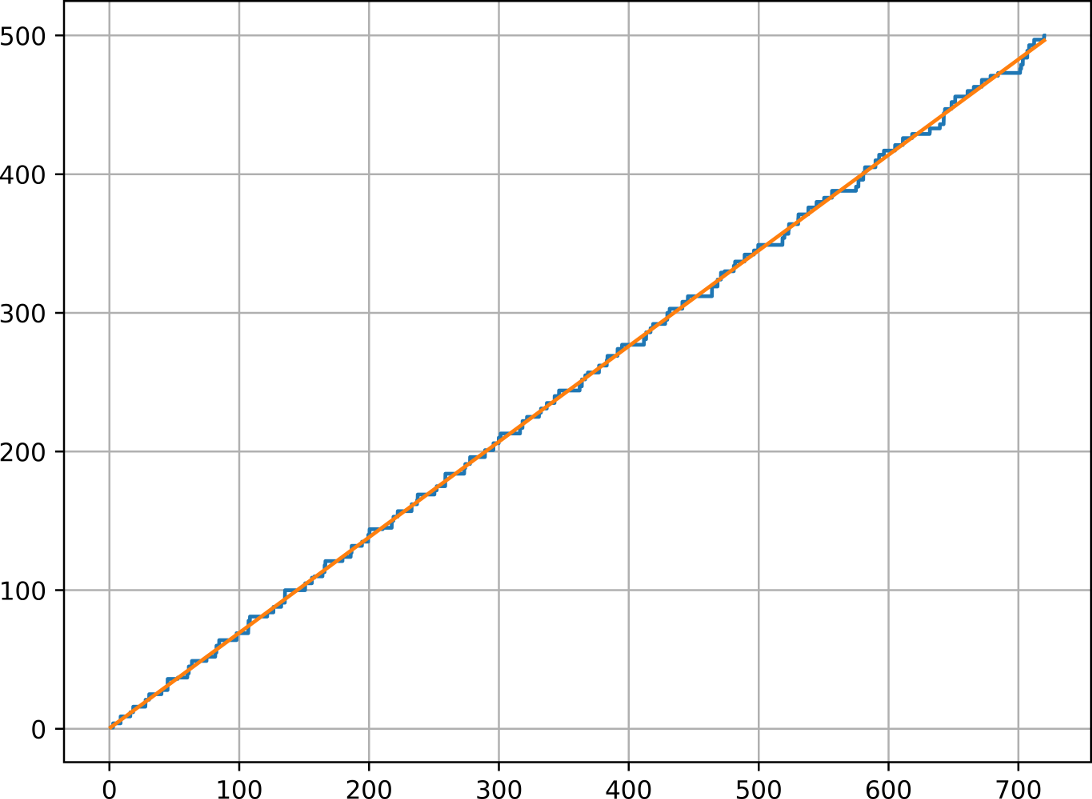}
    }
    \hfill
    \subfigure[The difference $D(t)$]
    {
        \includegraphics[width=2.1in]{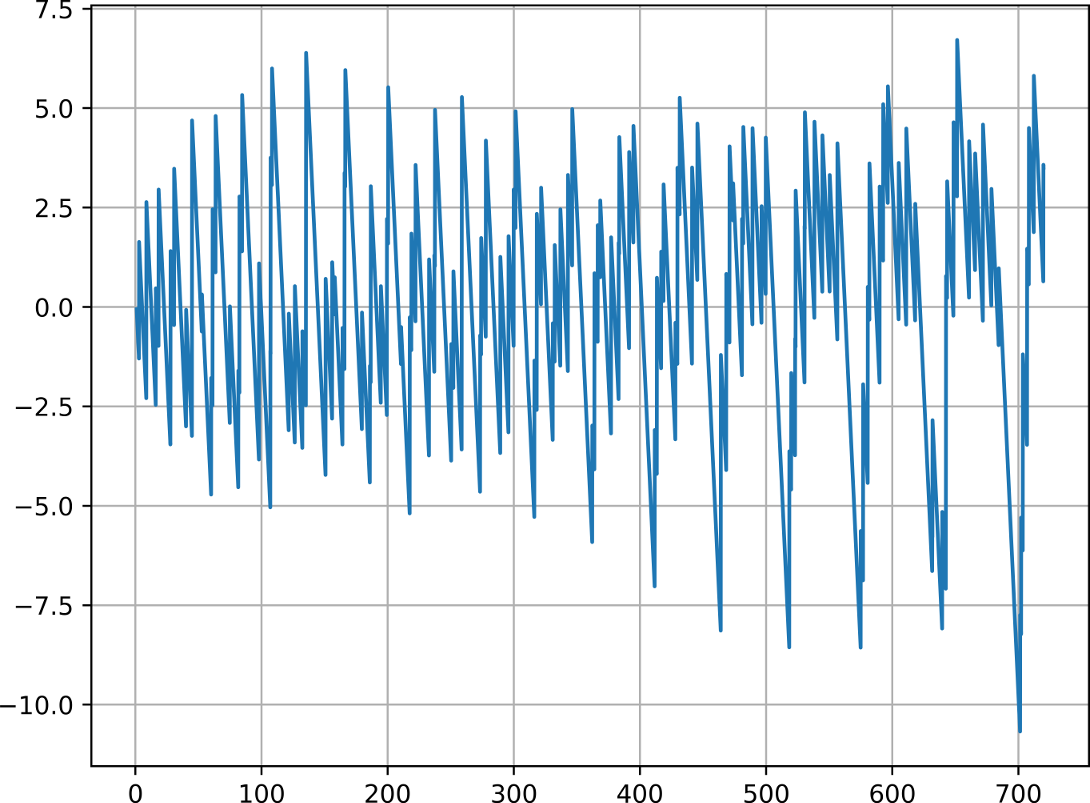}
    }
    \subfigure[The average of the difference $A(t)$]
    {
        \includegraphics[width=2.1in]{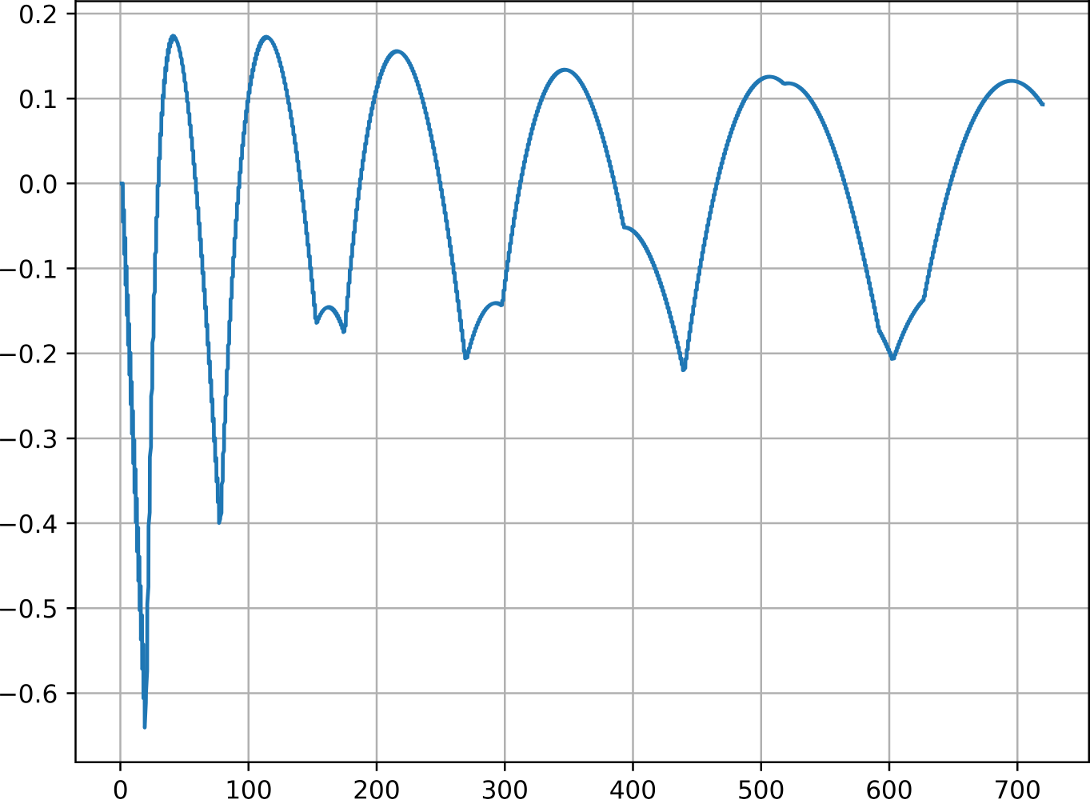}
    }
    \hfill
    \subfigure[The average after rescaling $g(t)$]
    {
        \includegraphics[width=2.1in]{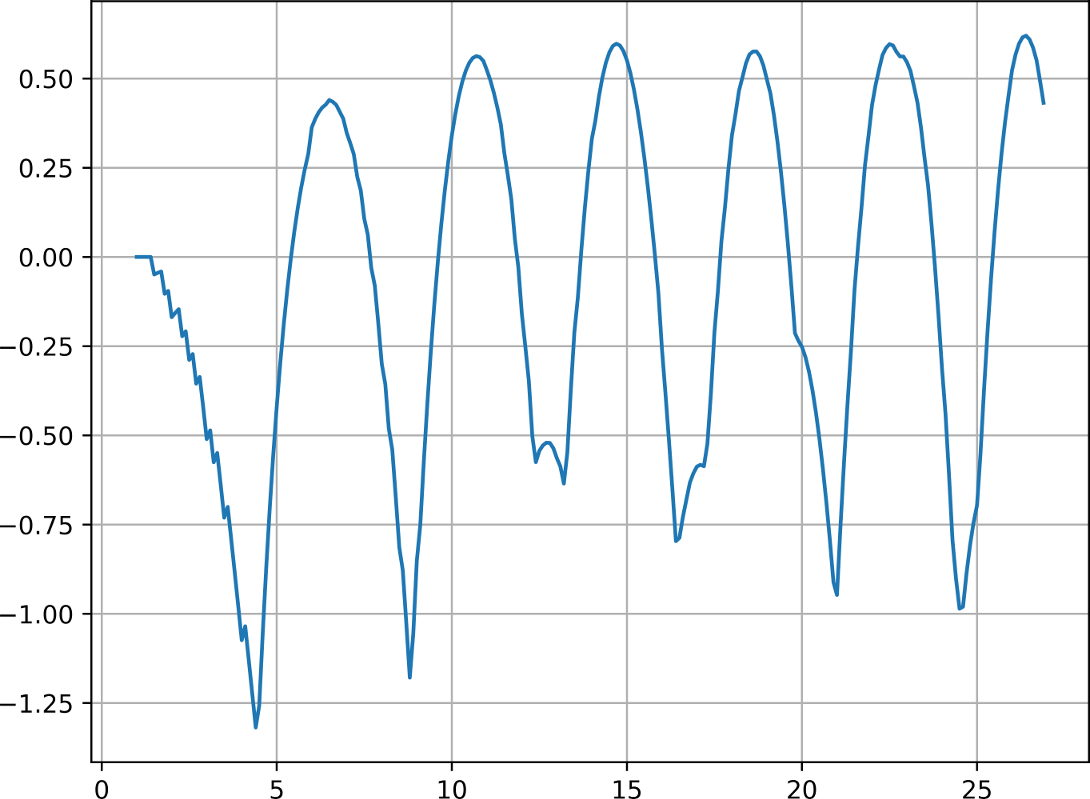}
    }
    \caption{Counting function on the icosahedron}
    \label{fig:icosgraphs}
\end{figure}

We again show the counting function $N(t)$, the remainder $D(t)$, and the averages  $A(t)$ and $g(t)$ in Figure \ref{fig:icosgraphs}. The constant $\frac{5\sqrt{3}}{4\pi}$ comes from the Weyl term and the constant $\frac{11}{30} = 12\cdot\frac{11}{360}$ is the value conjectured in [5].

\vspace{0.5cm}

\newpage

\section{The Cube.}

The cube has the same symmetry group as the octahedron (they are dual polyhedra), and the description of eigenfunctions for both polyhedra is similar. For each of the four $1$-dimensional representations we will find nonsingular eigenfunctions and we will be able to write down explicit formulas. Since the faces are squares rather than equilateral triangles, the formulas will be different.

\begin{figure}[h]
    \centering
    \subfigure[Straight]
    {
        \includegraphics[width=2in]{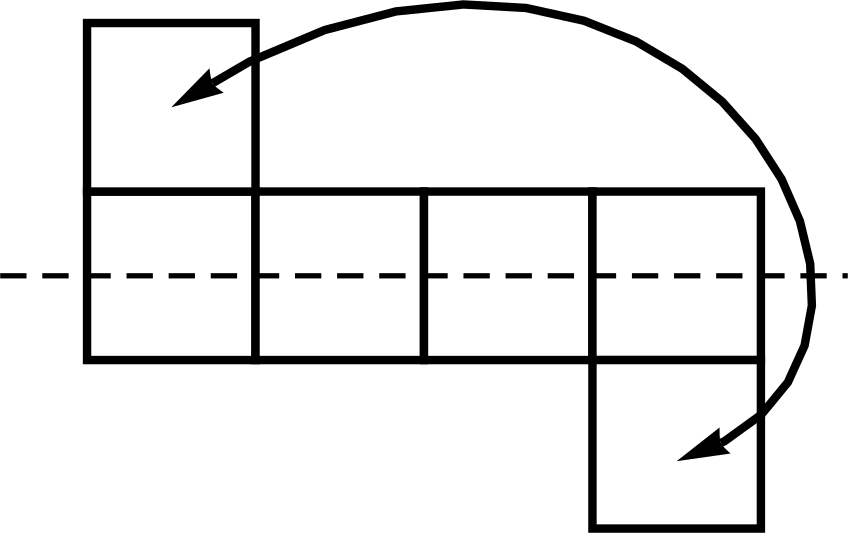}
    }
    \subfigure[Diagonal]
    {
        \includegraphics[width=2in]{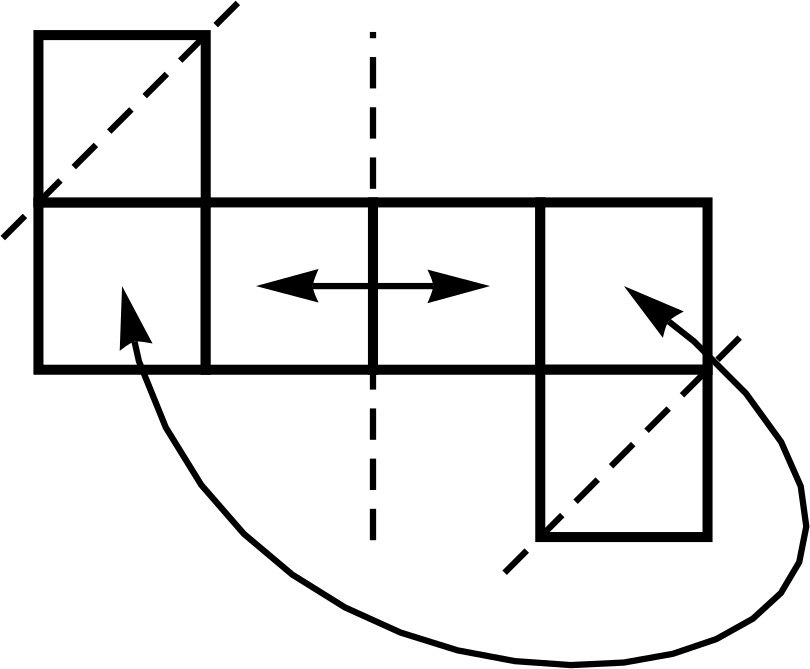}
    }
    \caption{Reflections on the cube}
    \label{fig:cuberefl}
\end{figure}

We distinguish two types of reflections, those that we call \emph{diagonal} that reflect in the diagonals of two opposite faces and switch the other faces in adjacent pairs, and those that we call \emph{straight} that reflect in horizontal or vertical bisectors of four faces and permute the remaining two opposite faces. These two types of reflections are illustrated in Figure \ref{fig:cuberefl}.

\begin{figure}[h]
    \centering
    \subfigure[  $1^{--}$]
    {
        \includegraphics[width=2.0in]{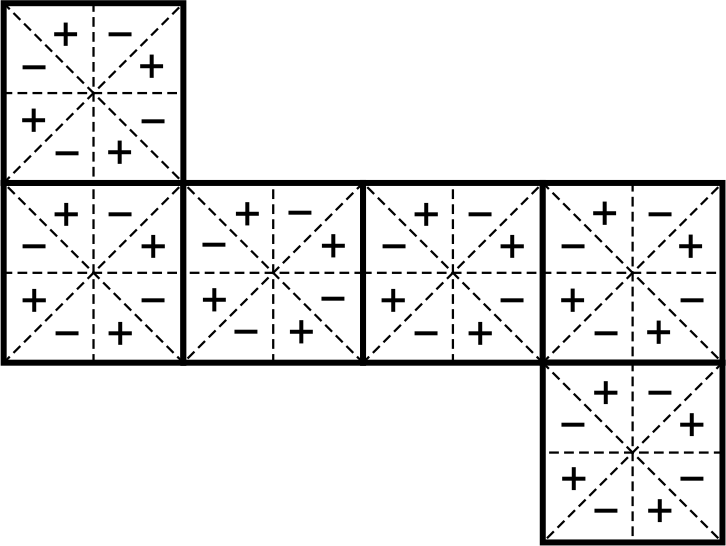}
    }
    \\
    \subfigure[$1^{+-}$ ]
    {
        \includegraphics[width=2.0in]{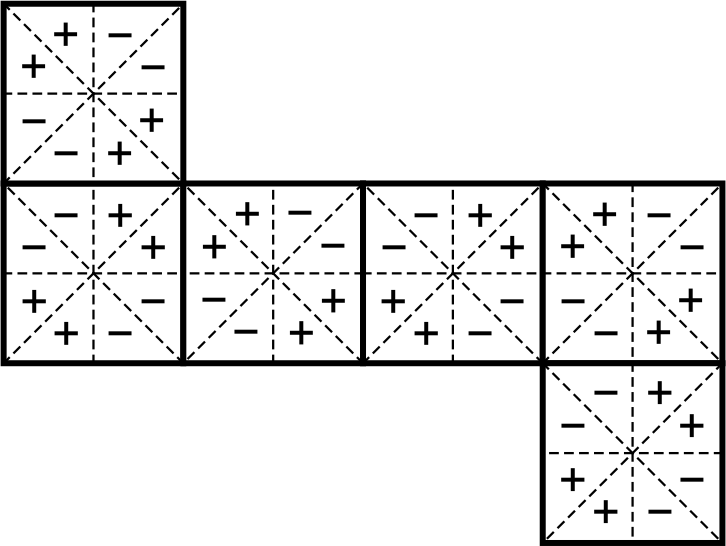}
    }
      \hspace{0.3in}
      \subfigure[$1^{-+}$ ]
    {
        \includegraphics[width=2.0in]{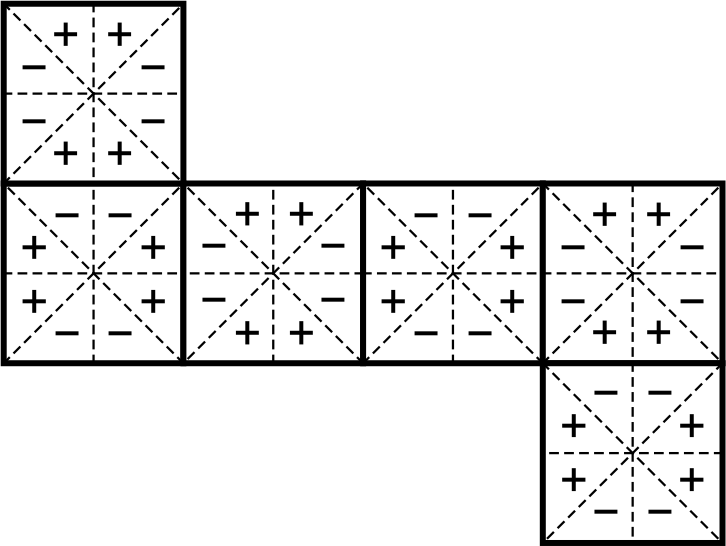}
    }
    \caption{Symmetric properties of eigenfunctions on the cube}
    \label{fig:cubesymm}
\end{figure}

We denote the corresponding symmetry types $1^{\pm\pm}$, where the first $\pm$ corresponds to symmetry $(+)$ and skew-symmetry ($-$) with respect to diagonal reflections, and the second $\pm$ refers to straight reflections. A function with any of these symmetry types is determined by its values on a single face, with extensions to the other faces obtained by even reflection for $1^{++}$ or $1^{+-}$, and odd reflection for $1^{-+}$ and $1^{--}$. We may choose as a fundamental domain a triangle one-eigth of a face. The reflection rules are illustrated in Figure \ref{fig:cubesymm}.

For simplicity we choose one face to be represented as the square in the x-y plane with vertices ($\pm \frac{1}{2},\pm \frac{1}{2})$. Since the half-turn $(x,y) \to (-x,-y)$ is the product of two reflections of the same type, all symmetry types produce eigenfunctions that are symmetric with respect to the half-turn. That means that all eigenfunctions are linear combinations of $\cos(\pi(kx+jy))$, for $(k,j)$ in the integer lattice $\mathbb{Z}^2$. Moreover, $k$ and $j$ must have the same parity. Indeed, for the $1^{++}$ and $1^{--}$ symmetry types we must have symmetry with respect to the rotations $x\to x+1$ and $y \to y+1$ (the product of two reflections, one in-face and one face-to-face), meaning $k$ and $j$ must both be even. For the $1^{+-}$ and $1^{-+}$ symmetry types we must have skew-symmetry with respect to these rotations, meaning $k$ and $j$ must both be odd. 

\newpage

\begin{figure}[h]
\centering
 \begin{minipage}{.45\textwidth}
 \centering
  \vspace{0.5cm}
  \includegraphics[width=1.55in]{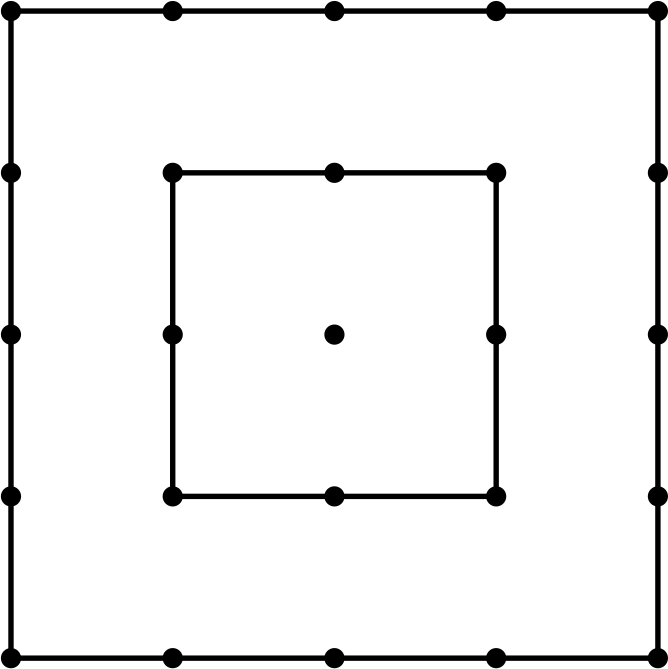}
  \vspace{0.6cm}
  \caption{Integer lattice}
  \label{fig:cubeintlat}
 \end{minipage}
 \begin{minipage}{.5\textwidth}
 \centering
  \includegraphics[width=2.6in]{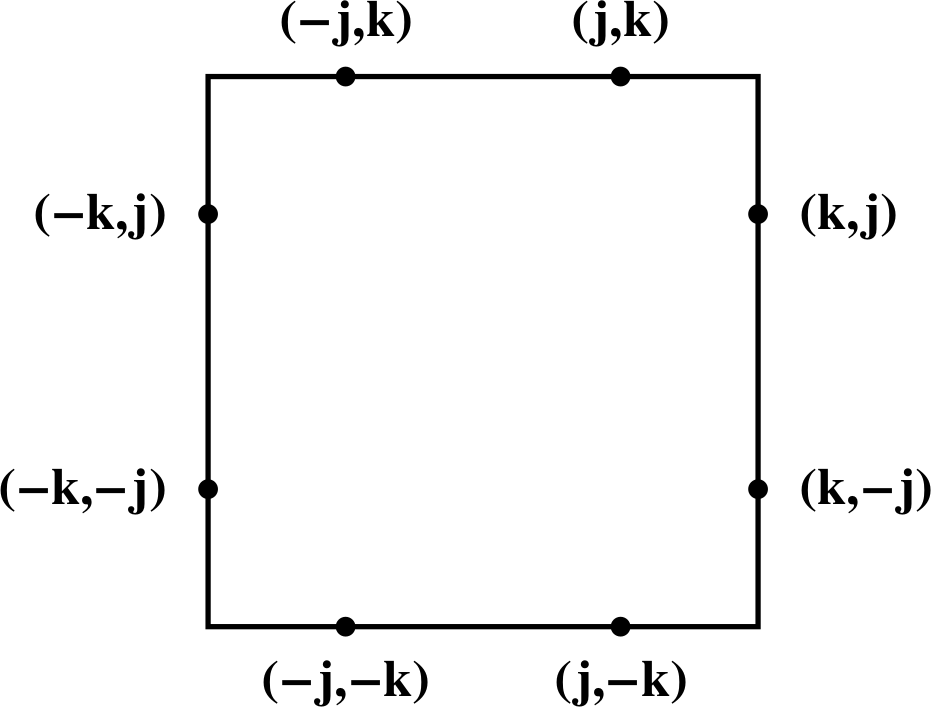}
  \caption{Generic orbit}
  \label{fig:cubegenorb}
 \end{minipage}
\end{figure}

To find the correct linear combinations we observe that the lattice decomposes into a disjoint union of concentric squares (aside from the origin), as shown in Figure \ref{fig:cubeintlat}.

The diagonal and straight reflections on the cube act on the lattice as diagonal and straight reflections. A generic orbit has eight points and is illustrated in Figure \ref{fig:cubegenorb}, with $k>j>0$ integers.

\newpage

\begin{figure}[h]
    \centering
    \subfigure[  $1^{--}$]
    {
        \includegraphics[width=1.33in]{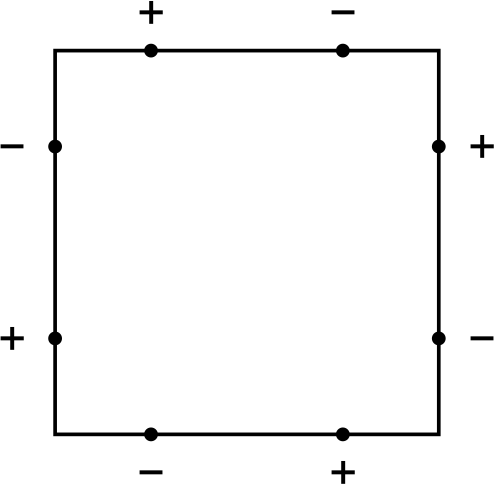}
    }
    \hspace{.15in}
    \subfigure[$1^{+-}$ ]
    {
        \includegraphics[width=1.33in]{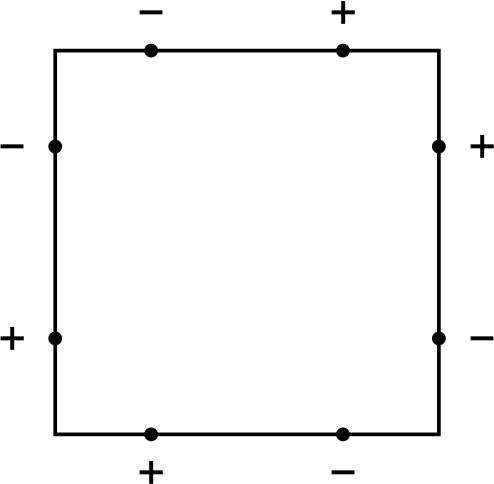}
    }
     \hspace{.15in}
      \subfigure[$1^{-+}$ ]
    {
        \includegraphics[width=1.33in]{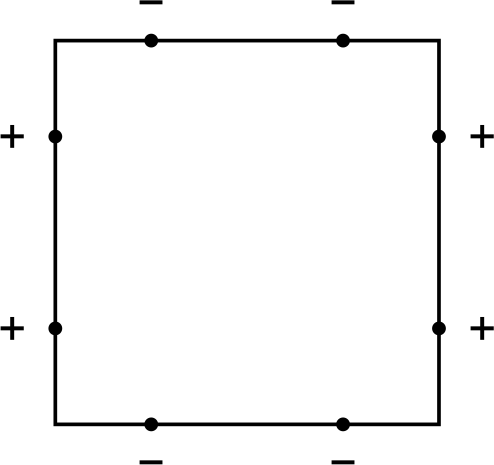}
    }
    \caption{Distribution of signs over the orbit on the cube}
    \label{fig:cubes}
\end{figure}

The distribution of signs over the orbit for the $1^{--}$,$1^{+-}$, and $1^{-+}$ eigenfunctions is shown in Figure \ref{fig:cubes}. Note that diametrically opposite lattice points always have the same sign, which is consistent with the fact that only cosines appear in the trigonometric polynomials.

From the diagrams we may read off the explicit formulas for the eigenfunctions with eigenvalue $\pi^2(j^2+k^2)$:

\begin{equation}
\label{eq:fiveone}
u_{++} = \cos(\pi(kx+jy)) + \cos(\pi(jx+ky)) + \cos(\pi(kx-jy)) + \cos(\pi(jx-ky))
\end{equation}
\begin{equation}
\label{eq:fivetwo}
u_{--} = \cos(\pi(kx+jy)) - \cos(\pi(jx+ky)) - \cos(\pi(kx-jy)) + \cos(\pi(jx-ky))
\end{equation}
\begin{equation}
\label{eq:fivethree}
u_{+-} = \cos(\pi(kx+jy)) + \cos(\pi(jx+ky)) - \cos(\pi(kx-jy)) - \cos(\pi(jx-ky))
\end{equation}
\begin{equation}
\label{eq:fivefour}
u_{-+} = \cos(\pi(kx+jy)) - \cos(\pi(jx+ky)) + \cos(\pi(kx-jy)) - \cos(\pi(jx-ky))
\end{equation}

\vspace{0.1cm}

\noindent where $k,j$ are both even in $(\ref{eq:fiveone})$, $(\ref{eq:fivetwo})$ and both odd in $(\ref{eq:fivethree})$, $(\ref{eq:fivefour})$.

There are two types of nongeneric orbits. When $j=0$  and $k$ is even the lattice points are at the midpoints of the edges. In this case only $u_{++}$ survives, and collapse to two terms:

\begin{equation}
\label{eq:fivefive}
u_{++} = \cos(\pi kx) + \cos(\pi ky)
\end{equation}

When $j=k$ the lattice points are at the corners of the square. In this case only $u_{++}$ or $u_{+-}$ survive, and again collapse into two terms:

\begin{equation}
\label{eq:fivesix}
u_{++} = \cos(\pi k(x+y)) + \cos(\pi k(x-y))
\end{equation}
\begin{equation}
\label{eq:fiveseven}
u_{+-} = \cos(\pi k(x+y)) - \cos(\pi k(x-y)) 
\end{equation}

\noindent for $k$ even in $(\ref{eq:fivesix})$ and $k$ odd in $(\ref{eq:fiveseven})$.

\newpage

In Table $4$ we list the eigenvalues normalized by dividing by $\pi^2$, with the integer values bolded. We extrapolate and will continue using the extrapolated data.

\begin{table}[h!]
\centering
\makebox[0pt][c]{\parbox{\textwidth}{
\begin{minipage}[b]{0.32\hsize}
\centering
\begin{tabular}{| c | c |}
\hline
\# & Eigenvalue \\ 	
\hline 
  \textbf{1} & \textbf{0} \\
  2 & 0.42105 \\ 
  3 & 0.42171 \\
  4 & 0.42197 \\
  5 & 1.16475 \\
  6 & 1.16502 \\
  7 & 1.16512 \\
  8 & 1.42522 \\
  9 & 1.43001 \\  
  \textbf{10} & \textbf{2.00027} \\
  11 & 2.59432 \\ 
  12 & 2.60125 \\
  13 & 2.60384 \\
  14 & 2.67862 \\
  15 & 2.67925 \\
  16 & 2.68175 \\
  17 & 3.81367 \\
  18 & 3.81781 \\
  19 & 3.81940 \\
  \textbf{20} & \textbf{4.00067} \\
\hline 
\end{tabular}
\end{minipage}
\hfill
\begin{minipage}[b]{0.32\hsize}
\centering
\begin{tabular}{| c | c |}
\hline
\# & Eigenvalue \\ 	
\hline 
  21 & 4.52692 \\
  22 & 4.52697 \\ 
  23 & 4.54599 \\
  24 & 4.61381 \\
  25 & 4.61602 \\
  26 & 5.65888 \\
  27 & 5.66338 \\
  28 & 5.66512 \\
  29 & 6.13609 \\
  30 & 6.15305 \\
  31 & 6.63945 \\
  32 & 6.65077 \\
  33 & 6.65518 \\
  34 & 7.00648 \\
  35 & 7.02039 \\
  36 & 7.02786 \\
  \textbf{37} & \textbf{8.00428} \\
  38 & 8.05707 \\
  39 & 8.07340 \\
  40 & 8.07945 \\   
\hline 
\end{tabular}
\end{minipage}
\hfill
\begin{minipage}[b]{0.32\hsize}
\centering
\begin{tabular}{| c | c |}
\hline
\# & Eigenvalue \\ 	
\hline 
  41 & 8.70184 \\
  42 & 8.70209 \\ 
  43 & 8.71324 \\
  44 & 9.41359 \\
  45 & 9.44725 \\
  46 & 9.70349 \\
  47 & 9.71256 \\
  48 & 9.71602 \\
  49 & 9.96909 \\
  \textbf{50} & \textbf{10.00591} \\
  51 & 11.02694 \\
  52 & 11.03827 \\
  53 & 11.04246 \\
  54 & 11.39163 \\
  55 & 11.39266 \\
  56 & 11.42616 \\
  57 & 11.95575 \\
  58 & 11.96738 \\
  59 & 12.69329 \\
  60 & 12.72420 \\
\hline 
\end{tabular}
\end{minipage}
}}
\vspace{0.5cm}
\caption{Normalized eigenvalues on the cube}
\end{table}

\newpage

\begin{figure}[h!]
    \centering
    \subfigure[Counting function $N(t)$ (blue) and $\left(\frac{3}{2\pi}t + \frac{7}{18}\right)$ (orange)]
    {
        \includegraphics[width=2.1in]{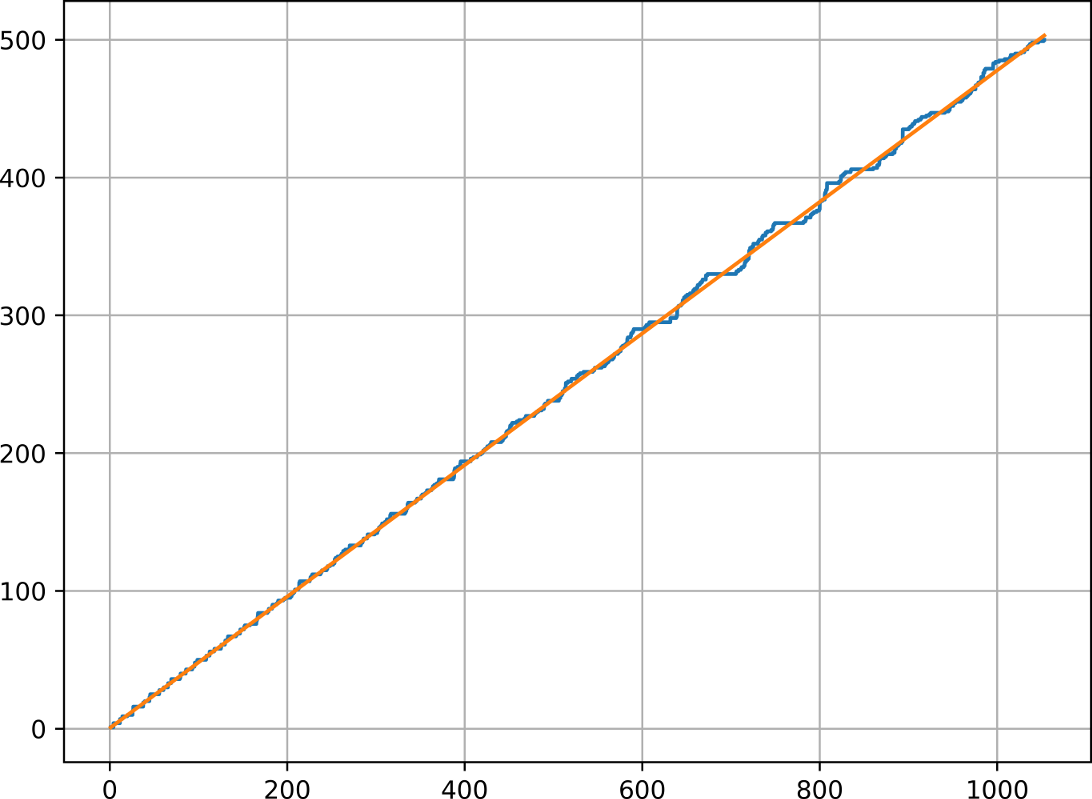}
    }
    \hfill
    \subfigure[Remainder $D(t)$]
    {
        \includegraphics[width=2.1in]{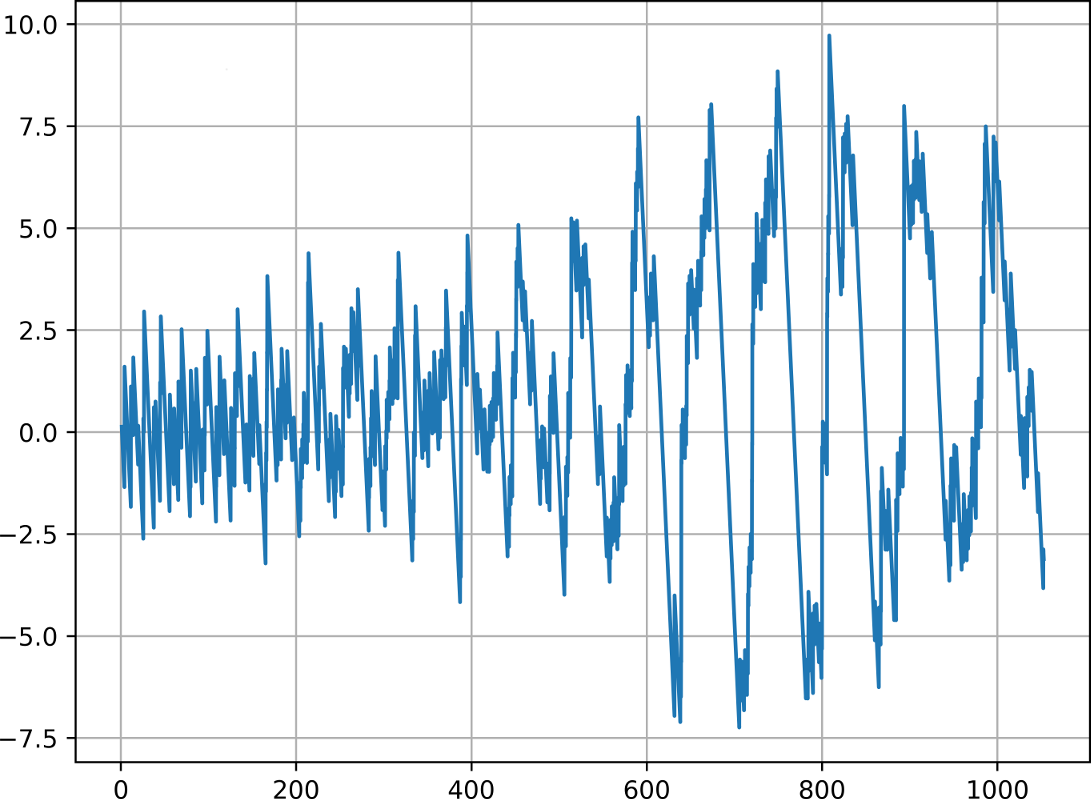}
    }
    \subfigure[The average of the difference $A(t)$]
    {
        \includegraphics[width=2.1in]{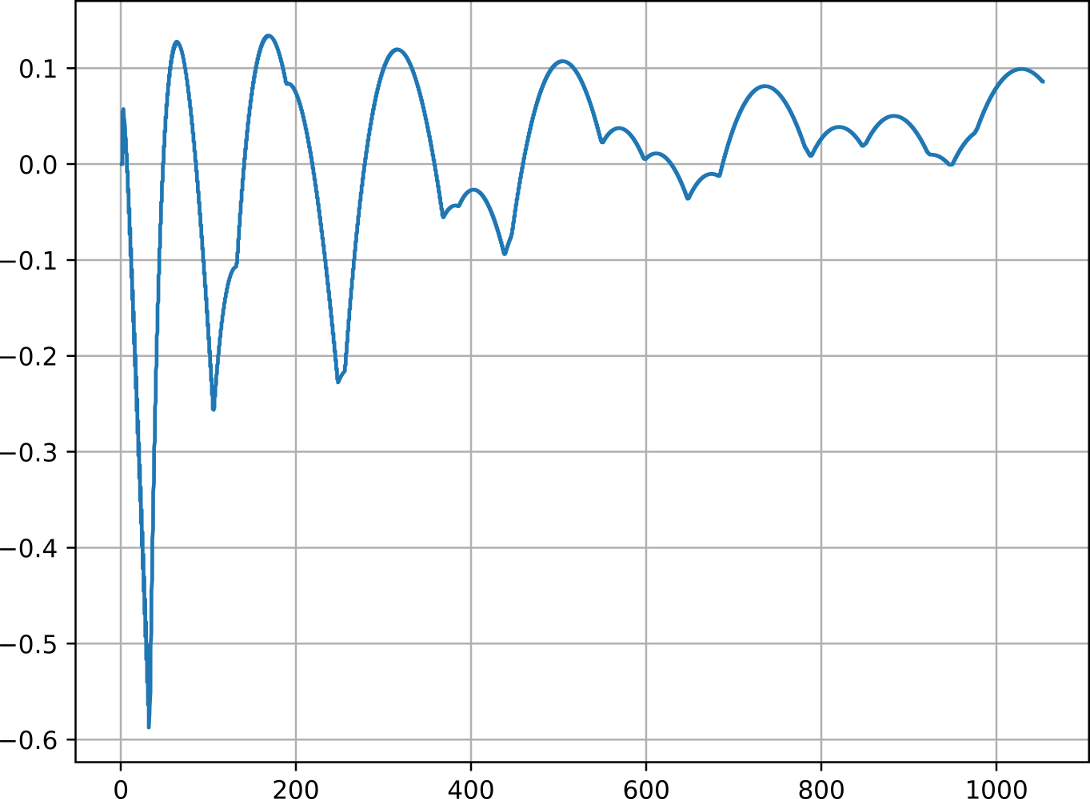}
    }
    \hfill
    \subfigure[The average after rescaling $g(t)$]
    {
        \includegraphics[width=2.1in]{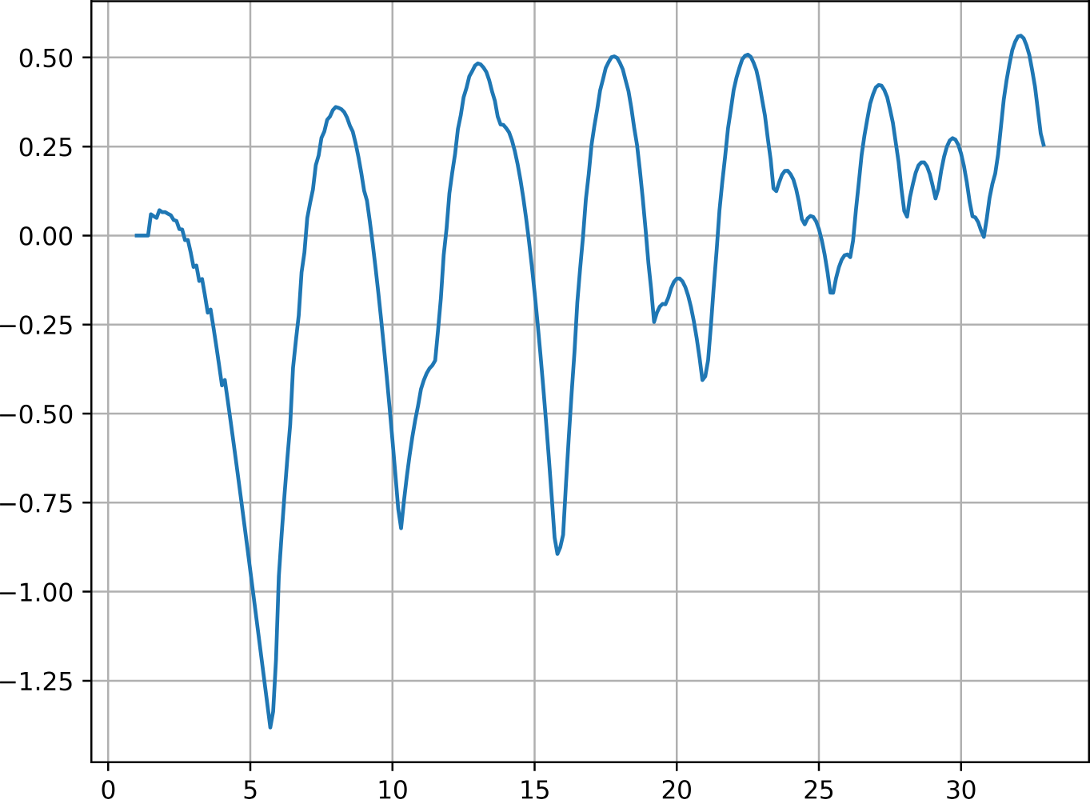}
    }
    \caption{Counting function on the cube}
    \label{fig:cubegraphs}
\end{figure}

We again show the counting function $N(t)$, the remainder $D(t)$, and the averages $A(t)$ and $g(t)$ in Figure \ref{fig:cubegraphs}. The constant $\frac{3}{2\pi}$ comes from the Weyl term and the constant $\frac{7}{18}=8\cdot\frac{7}{144}$ is the value conjectured in [5].

\vspace{0.4cm}

\clearpage


\bibliographystyle{amsplain}

\begin{thebibliography}{10}

\bibitem {JS} S. Jayakar, R. S. Strichartz, \textit{Average number of lattice points in a disk}, Commun. Pure Appl. Anal., \textbf{15} (2016), 1--8

\bibitem {Sen} A. N. Sengupta, \textit{Representing Finite Groups: A semisimple introduction}, Springer, New York, 2012

\bibitem {Se} J. Serre, \textit{Linear Representations of Finite Groups}, Springer, New York, 1977

\bibitem {Sh} T. Shioya, \textit{Geometric Analysis on Alexandrov Spaces}, Sugaku Expositions, \textbf{24} (2011), 145--167

\bibitem {Str} R. S. Strichartz, \textit{Average error for spectral asymptotics on surfaces}, Commun. Pure Appl. Anal., \textbf{15} (2016), 9--39

\bibitem {W} R. S. Strichartz, S. C. Wiese, \textit{Spectrum of the Laplacian on regular polyhedra}, \url{http://pi.math.cornell.edu/\~polyhedral} (updated August 7, 2018)

\end{thebibliography}

\end{document}